\documentclass{conm-p-l}
\usepackage{array,epsfig,latexsym,tabularx}

\begin{document}


\newtheorem{example}{Example}[section]
\newtheorem{note}[example]{Note}
\newtheorem{theorem}[example]{Theorem}
\newtheorem{corollary}[example]{Corollary}
\newtheorem{definition}[example]{Definition}
\newtheorem{proposition}[example]{Proposition}
\newtheorem{algorithm}[example]{Algorithm}
\newtheorem{lemma}[example]{Lemma}
\newtheorem{problem}[example]{Problem}
\newtheorem{conjecture}[example]{Conjecture}


\newcommand{\rai}{\rightarrow \, \infty}
\newcommand{\infinity}{\infty}
\renewcommand{\mod}{\mbox{mod}\,}

 
\font\twelvesym=msbm10 at 12pt
\font\tensym=msbm10
\font\sevensym=msbm7
\font\fivesym=msbm5
\newfam\ssymfam
\textfont\ssymfam=\tensym
\scriptfont\ssymfam=\sevensym
\scriptscriptfont\ssymfam=\fivesym
\def\ssym{\fam\ssymfam\tensym}


\newcommand{\Z}{{\ssym Z}}
\newcommand{\N}{{\ssym N}}

\newcommand{\B}{{\mathcal B}}
\newcommand{\D}{{\mathcal D}}
\renewcommand{\H}{{\mathcal H}}
\newcommand{\M}{{\mathcal M}}
\renewcommand{\P}{{\mathcal P}}
\newcommand{\V}{{\mathcal V}}
\newcommand{\Sc}{{\mathcal S}}

\newcommand{\wt}{{\rm wt\,}}
\newcommand{\bwt}{wt \,}
\newcommand{\Proof}{\medskip\noindent {\it Proof: }}
\newcommand{\cqfd}{\hfill $\Box$ \medskip}
\newcommand{\boldm}{\mbox{\boldmath$m$}}
\newcommand{\boldn}{\mbox{\boldmath$n$}}
\newcommand{\bolde}{\mbox{\boldmath$e$}}
\newcommand{\boldu}{\mbox{\boldmath$u$}}
\newcommand{\boldQ}{\mbox{\boldmath$Q$}}
\newcommand{\boldR}{\mbox{\boldmath$R$}}
\newcommand{\boldC}{\mbox{\boldmath$C$}}
\newcommand{\sboldm}{\mbox{\boldmath$\scriptstyle m$}}
\newcommand{\sboldn}{\mbox{\boldmath$\scriptstyle n$}}
\newcommand{\sbolde}{\mbox{\boldmath$\scriptstyle e$}}
\newcommand{\sboldu}{\mbox{\boldmath$\scriptstyle u$}}
\newcommand{\sboldQ}{\mbox{\boldmath$\scriptstyle Q$}}
\newcommand{\sboldR}{\mbox{\boldmath$\scriptstyle R$}}
\newcommand{\sboldC}{\mbox{\boldmath$\scriptstyle C$}}
\newcommand{\boldlambda}{\mbox{\boldmath$\lambda$}}
\newcommand{\wombat}{\rule[-6pt]{0pt}{46pt}}






\def\makestrut#1#2{{\dimen12=#2
\divide\dimen12 by 4\dimen11=\dimen12\multiply\dimen11 by 3
\global\setbox#1=\hbox{\vrule height\dimen11 depth\dimen12 width0pt}}}

\newdimen\tadhdimen \newdimen\tabhdimen \newdimen\vdimen
\newdimen\smtadhdimen \newdimen\smtabhdimen
\newcount\splurt
\newbox\tadstrut \newbox\tabstrut
\newbox\smtadstrut \newbox\smtabstrut


\def\setyoungsize#1#2{          
  \tadhdimen=#1\tabhdimen=#1\advance\tabhdimen by -0.4truept%
  \vdimen=#2%
  \makestrut\tadstrut\vdimen
  \advance\vdimen by -0.4pt%
  \makestrut\tabstrut\vdimen}


\setyoungsize{18pt}{13pt}


\def\youngt#1{%
  \vcenter{\offinterlineskip
  \halign{&\copy\tadstrut\hbox to \tadhdimen{\hss$##$\hss}\cr #1}}}
\def\youngd#1{%
  \vcenter{\offinterlineskip
  \halign{&\vrule##&\copy\tabstrut\hbox to \tabhdimen{\hss$##$\hss}\cr #1}}}


\hyphenation{boson-ic 
             ferm-ion-ic 
	     para-ferm-ion-ic
             two-dim-ension-al
	     two-dim-ension-al}

\title{Melzer's identities revisited}
\thanks{Research supported by the Australian Research Council (ARC)}
\dedicatory{Dedicated to Professor George Andrews on the occasion
            of his 60th birthday.}

\author{Omar~Foda}
\address{Department of Mathematics and Statistics,
         The University of Melbourne,
         Park\-ville, Victoria 3052, Australia.}
        \email{foda@maths.mu.oz.au} 
	
\author{Trevor~A.~Welsh} 
\address{Department of Mathematics and Statistics, 
         The University of Melbourne, 
         Park\-ville, Victoria 3052, Australia.} 
	\email{trevor@maths.mu.oz.au} 

\begin{abstract}

We further develop the finite length path generating transforms 
introduced previously, and use them to obtain constant sign 
polynomial expressions that reduce, in the limit of infinite path 
lengths, to parafermion and ABF Virasoro characters. This provides
us, in the ABF case, with combinatorial proofs of Melzer's 
boson-fermion polynomial identities.

\end{abstract}

\maketitle

\newpage

\setcounter{secnumdepth}{10}

\section{Introduction}

One can think of exactly solvable models, in statistical mechanics 
\cite{baxter-book} and in quantum field theory \cite{dms-book}, as 
concrete realisations of certain mathematical structures. These 
structures are so powerful, that they allow us to compute, at least 
in principle, an infinite number of physical quantities in each 
solvable model. Computing {\sl one} such quantity suffices to call 
the corresponding model {\sl solved}. 

Of particular interest are the connections between exact solutions 
and infinite dimensional algebras \cite{jimbo-miwa-book}. One aspect
of this connection is the observation, first made in \cite{djkmo}, 
that the one-point functions of regime-III restricted solid-on-solid 
ABF models \cite{abf}, suitably normalised, turn out to be characters 
of Virasoro highest weight modules \cite{bpz}\footnote{For the rest 
of this work, we refer to the regime-III restricted solid-on-solid ABF 
models simply as \lq ABF models\rq. The spectrum generating algebra of the 
ABF models is the Virasoro algebra of \cite{fqs}, with central charge 
$c=1-6/p(p+1), p=3, 4, \cdots$}. 

Similarly, the one-point functions of regime-II ABF models are characters 
of parafermion highest weight modules\footnote{For the rest of this work, 
we refer to regime-II restricted solid-on-solid ABF models simply as 
\lq parafermion models\rq. The spectrum generating algebra of the
parafermion models is the $\Z$ algebra of \cite{lepowsky-wilson},
with central charge $c=2(p-2)/(p+1), p=3, 4, \cdots$}. 

The purpose of this work is to discuss combinatorial aspects of the ABF 
and parafermion one-point functions, or equivalently the corresponding 
Virasoro and parafermion characters\footnote{As we will see below, the
ABF and parafermion models are related by the transformation
$q \rightarrow q^{-1}$, where $q$ is the nome of the elliptic functions
used to parametrise the two-dimensional weights of these models, or 
equivalently, the expansion parameter that appears in the $q$-series
expression of the characters. For that reason, for each statement that 
we make about the ABF models, a corresponding statement can be made 
about the parafermion models.}.
The shortest route to the combinatorics that we are interested in is 
{\it via} the statistical mechanical side of the problem.

{}From the lattice point of view, a one-point function is the normalised 
generating function of an infinite set of two-dimensional configurations 
with very complicated weights (typically products of trigonometric or even 
elliptic functions). Baxter's corner transfer matrix method reduces the 
above problem to computing the generating function of an infinite set of 
one-dimensional configurations with relatively very simple weights (typically 
simple powers of a parameter $q$) \cite{baxter-book}.

{}From a combinatorial point of view, the set of weighted one-dimensional 
configurations is the starting point of this work. One does not need to 
know anything about the underlying physical models, or their connections
with infinite dimensional algebras\footnote{We refer the reader to 
\cite{jimbo-miwa-book} for an excellent introduction to the algebraic 
approach to exactly solvable lattice models.}.
We are handed a set of one-dimensional combinatorial objects, rules for 
computing their weights, and the task of computing their generating 
functions\footnote{We refer to computing 
the generating function of a set of weighted object simply as 
{\sl $q$-counting}.}.

There is no unique method to compute the generating function of 
a weighted set. Different methods produce different expressions. 
Since they all represent the same generating function, equating 
them produces $q$-series identities. 

One way to compute a generating function is \lq sieving\rq, or 
inclusion-exclusion \cite{andrews-red-book}. By construction, this 
method produces an expression whose terms have alternating signs.
In other words, the coefficients of the $q$-series so expressed,
are not manifestly positive definite.
However, we know that they are positive definite, since we are counting 
objects. The alternating-sign $q$-series expressions for the
one-dimensional configurations coincide with the Rocha-Caridi
expressions for the Virasoro characters \cite{rocha-caridi}. 

In the context of the Virasoro characters, the Stony Brook group 
were the first to conjecture that there exist constant-sign $q$-series 
expressions\footnote{For a complete listing of the papers of the Stony 
Brook group on this subject, we refer the reader to 
\cite{stony-brook-review}.}. For physical reasons that are beyond 
the scope of this work, the alternating-sign expressions are referred 
to as {\sl bosonic}. The constant-sign expressions are referred to as 
{\sl fermionic}. On equating these expressions, one obtains 
{\sl boson-fermion $q$-series identities}.

One approach to proving such identities, is to work at the level of 
finite versions of the combinatorial objects under consideration.
$q$-Counting these finite sets produces boson-fermion $q$-polynomial 
identities. Since the initial conjectures of the Stony Brook group, 
there has been many further conjectures and proofs of $q$-polynomial 
identities. In this work, we are interested in Melzer's polynomial
identities \cite{melzer}. For each one-point function, of each ABF 
model, Melzer conjectured four boson-fermion $q$-polynomial 
identities\footnote{For each one-point function, these identities
are not independent. As we will see below, two of them are related 
to another two, that belong to another one-point function, by means
of a simple up-down reflection of the combinatoric objects that 
are counted.}.
Because of their relative simplicity, these conjectures have served 
as ideal testing grounds for various approaches towards proving 
boson-fermion identities. Proofs of a subset of Melzer's identities 
were obtained in \cite{berkovich,berkovich-mccoy,bms} using recursion 
techniques. A complete proof using the same methods is given in 
\cite{schilling}. Combinatorial proofs of a subset of these 
identities were obtained in 
\cite{dasmahapatra-foda, foda-warnaar, warnaar1, warnaar2}.

In this work, we are interested in a combinatorial proof of 
the full set of Melzer's identities. We obtain such a proof 
by extending our previous work on path generating transforms
\cite{flpw}. Though, strictly speaking, we do not obtain new 
final results, the method that we use is new. We hope that 
this method gives further insight into the combinatorics of 
Virasoro highest weight modules, which have turned out to be 
such rich and fascinating objects.

\subsection{Outline of paper}

In Section \ref{PathDefSec}, we define the combinatorial objects,
called {\sl paths}, that we are interested in $q$-counting. In 
sections \ref{ABFPathSec} and \ref{ParaPathSec}, we define, in 
terms of paths, the two generating functions of most importance 
to us\footnote{In \cite{abf}, these generating functions are 
denoted $X_L(a,b,c)$. and $x_L(a,b,c)$.}.

The first, $X^{p'}_L(a,b,c)$, gives the finitised one-point function 
of the ABF models. The second, $x^{p'}_L(a,b,c)$, gives the finitised 
one-point function of the parafermion models. Lemma \ref{DualityLem} 
states the relationship between $X^{p'}_L(a,b,c)$ and 
$x^{p'}_L(a,b,c)$\footnote{In \cite {abf}, $x^{p'}_L(a,b,c)$ 
is defined through its relationship with $X^{p'}_L(a,b,c)$, and not 
directly in terms of the paths as we do here.}. 

However, instead of the two functions $X^{p'}_L(a,b,c)$ and 
$x^{p'}_L(a,b,c)$, we prefer to work with certain renormalisations 
thereof, which we denote $\chi^{p'-1, p'}_{a,b,c}(L)$ and
$\chi^{1, p'}_{a,b,c}(L)$ 
respectively\footnote{Here, we retain the notation of \cite{flpw}.}.
We make this change for two reasons. Firstly, the analysis here then 
parallels that of \cite{flpw} and a comparison can be readily made. 
Moreover, the techniques of this work and those of \cite{flpw} are 
then ready to combine so that the other cases of \cite{forrester-baxter} 
may be investigated. Secondly, the connection with partitions satisfying
prescribed hook-difference conditions, analysed in \cite{abbbfv}, is 
then apparent. We discuss this in Appendix B, relying somewhat on 
the analysis of \cite{flpw}.

Section \ref{AltPresSec} indicates how $\chi^{1,p'}_{a,b,c}(L)$ may be
directly determined from the paths. Our strategy is to obtain 
expressions for $\chi^{1,p'}_{a,b,c}(L)$ using combinatorial techniques 
applied to the paths, and then to obtain expressions for the other 
functions mentioned above from these. We choose to work with these 
parafermion models for compatibility with our work in \cite{flpw}.
Furthermore, combining the techniques with those of \cite{flpw} will, 
in future work,
enable further models from \cite{forrester-baxter} to be tackled.

In Section \ref{AltPresSec}, each vertex of a path is designated 
either {\sl scoring} or {\sl non-scoring}. In Section \ref{StrikingSec}, 
the {\sl striking sequence}\footnote{Although a similar notion, the 
definition of the striking sequence given in \cite{flpw} differs from 
that used here.} of a path is defined, and the means to designate the 
first point and last point of a path as scoring or non-scoring is 
given. In Section \ref{RestGenSec}, we define the generating function 
$\chi^{1,p'}_{a,b,e,f}(L,m)$ for paths having a certain length, a certain 
number of non-scoring vertices, and first and last vertices of a 
certain nature.

In Section \ref{TranSec}, we introduce the cornerstone of our method.
This is the notion of a transform which enables us to express the 
generating functions in terms of those of a \lq simpler\rq\ model.
This transform is called a ${\B}$-transform.

Following the action of a ${\B}$-transform, the path may be 
extended by adding a number of segments, alternating in direction, to 
the left end. This process, which we refer to (for physical reasons 
that do not concern us here) as {\sl inserting particles}, is described 
in Section \ref{InsertSec}. These particles are then allowed to 
{\em move} through the path, as described in Section 
\ref{MovesSec}\footnote{Our ${\B}$-transform, and the subsequent 
insertion of particles are all direct extensions of ideas that we 
learnt from \cite{agarwal-bressoud, bressoud}. The \lq particle moves\rq\ 
also appear in \cite{agarwal-bressoud, bressoud} in the context of 
somewhat different models. They appear in the context of the ABF 
models in \cite{warnaar1, warnaar2}.}. As shown in Section \ref{WavesSec}, 
this whole process enables $\chi^{1,p'}_{a,b,e,f}(L,m)$ to be expressed 
in terms of various $\chi^{1,p'-1}_{a',b',e,f}(m,m')$.

Using the techniques of the previous sections, constant-sign expressions 
for $\chi^{1,p'}_{a,b,c}(L)$ are obtained in Section \ref{ParaFunSec} by 
employing a succession of ${\B}$-transforms. It turns out that this may 
be accomplished in four different ways and these lead to four different 
constant-sign expressions. Two of these are derived in Sections 
\ref{FirstSystemSec} and \ref{SecondSystemSec}. The two constant-sign 
expressions that result are necessarily equal, although this is by no means 
obvious\footnote{An analytic proof of this fact is obtained in 
\cite{schilling}.}. Using a simple symmetry argument, these expressions 
yield the other two constant-sign expressions.

Finally, in Section \ref{ABFFunSec}, the relationship between
the ABF and parafermion models is employed to obtain the 
constant-sign expressions for the ABF one-point functions,
that had been conjectured by Melzer \cite{melzer}.

In Appendix A, we describe the \boldm\boldn-systems that pertain 
to the constant-sign expressions obtained. In Appendix B, we discuss 
the aforementioned connection with partitions satisfying prescribed 
hook-difference conditions.

\section{Combinatorics of highest weight modules}

\subsection{Paths}\label{PathDefSec}

Let $p'\in\N$ with $p'\ge2$.
A path $h$ of length $L$ is a sequence $h_0,h_1,h_2,\ldots,h_L,$
of integer heights such that $1\le h_i\le p'-1$ for
$0\le i\le L$, and such that $h_{i+1}=h_i\pm1$ for $0\le i<L$.
Such paths may readily be depicted on a two-dimensional
$L\times (p'-1)$ grid. The path is then the series of contiguous
line segments passing from $(i,h_i)$ to $(i+1,h_{i+1})$ for
$0\le i<L$. Note that each of these line segments is either
in the NE direction or in the SE direction.
It will be useful to define the length function: $L(h)=L$.

The following is a typical path $h$.  It's length is $L(h)=11$.

\medskip
\centerline{\epsfig{file=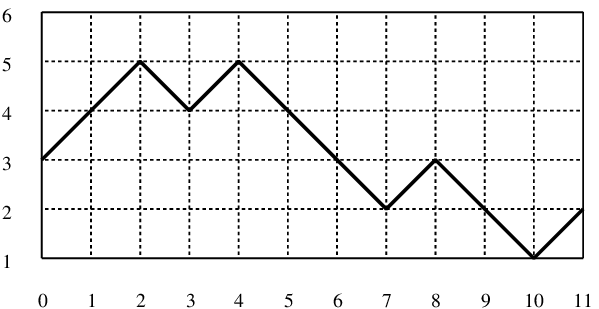}}
\medskip
\centerline{Figure\ 1.}
\medskip

\noindent
If $L+a-b$ is even, define ${\P}^{p'}_{a,b}(L)$ to be the set
of all paths $h$ of length $L$ with $h_0=a$ and $h_L=b$.

In \cite{abf} and \cite{forrester-baxter}, a number of ways of assigning 
a weight to each path are described. In this paper, we are interested in 
only two of these. They are the cases considered in \cite{abf}, and therein 
denoted regime III and regime II\footnote{In \cite{forrester-baxter}, the 
many ways of assigning weights are indexed by $p$ with $1\le p<p'$ and $p$ 
coprime to $p'$ (see also \cite{flpw}). Regime III of \cite{abf} is then 
the case $p=p'-1$ and regime II of \cite{abf} is the case $p=1$.}. As we 
mentioned above, we shall refer to these as the ABF model and the parafermion 
model respectively.

In each case, a weight is assigned to each path $h$ only
after an {\it extra} point $h_{L+1}$ satisfying
$1\le h_{L+1}\le p'-1$ and $h_{L+1}=h_L\pm1$ is specified.
Then, if $1\le a,b,c\le p'-1$ with $c=b\pm1$ and $L\ge0$ is such that
$L+a-b$ is even, we define ${\P}^{p'}_{a,b,c}(L)$ to be the set%
\footnote{We maintain a clear distinction between
${\P}^{p'}_{a,b}(L)$ and ${\P}^{p'}_{a,b,c}(L)$.
Namely, $h_{L+1}$ is defined for each element $h$ of the latter set,
whereas it is not for the former set.
This implies that the $L$th vertex of $h$ has a definite shape for each
element of the latter, but not for the former.}
of all paths $h$ of length L, such that $h_0=a$, $h_L=b$ and $h_{L+1}=c$.

\subsection{ABF}\label{ABFPathSec}

In regime III of \cite{abf} (the ABF model),
each path $h$ is assigned a weight $\bwt^{\rm III}(h)$ given by:
\begin{equation}\label{BwtIIIDef}
\bwt^{\rm III}(h)=\sum_{i=1}^L ic^{\rm III}(h_{i-1},h_{i},h_{i+1}),
\end{equation}
where, the function $c^{\rm III}(h_{i-1},h_{i},h_{i+1})$
is defined by:
\begin{eqnarray*}
c^{\rm III}(h-1,h,h+1)&=&1/2\,;\\
c^{\rm III}(h+1,h,h-1)&=&1/2\,;\\
c^{\rm III}(h-1,h,h-1)&=&0\,;\\
c^{\rm III}(h+1,h,h+1)&=&0\,.
\end{eqnarray*}
Note that these four cases correspond to the four different
vertex shapes. They appear as follows.

\medskip
\centerline{\epsfig{file=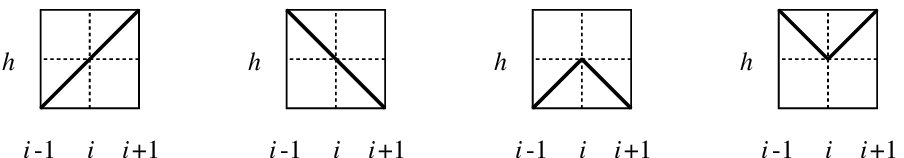}}
\medskip
\centerline{Figure\ 2.}
\medskip
 
\noindent
They will be referred to as a {\em straight-up} vertex,
a {\em straight-down} vertex, a {\em peak-up} vertex and
a {\em peak-down} vertex respectively.

The generating function for paths in the ABF model is defined to be:
\begin{equation}\label{GenIIIDef}
X^{p'}_L(a,b,c;q)=\sum_{h\in{\P}^{p'}_{a,b,c}(L)} q^{\bwt^{\rm III}(h)}.
\end{equation}
We set $X^{p'}_L(a,b,c)=X^{p'}_L(a,b,c;q)$.
A bosonic (i.e. alternating sign) expression for
$X^{p'}_L(a,b,c)$ is obtained in \cite{abf} (Theorem 2.3.1.).
It gives
\begin{equation}\label{RenormABFEq}
X^{p'}_L(a,b,c)=q^{-{1\over4}(a-b)(a-c)}\chi^{p'-1,p'}_{a,b,c}(L),
\end{equation}
where:
\begin{eqnarray}\label{FinRochaEq}
\chi^{p'-1,p'}_{a,b,c}(L)&=&
\sum_{\lambda=-\infinity}^\infinity
q^{\lambda\left((p'-1)(\lambda p'-a)+p'r\right)}
\left[ {L\atop {L+a-b\over2}-p'\lambda} \right]_q\nonumber\\[0.5mm]
&&\qquad\quad
-\sum_{\lambda=-\infinity}^\infinity
q^{(\lambda p'-\lambda+r)(\lambda p'+a)}
\left[ {L\atop {L+a-b\over2}-p'\lambda-a} \right]_q,
\end{eqnarray}
where $r=(b+c-1)/2$ (i.e.\ $r=\min(b,c)$), and
as usual, the Gaussian polynomial
$\left[ {A \atop B} \right]_q$ is defined to be:
\begin{equation}\label{gaussian}
\left[ {A \atop B} \right]_q =
     \frac{\prod_{i=1}^{A  }(1 - q^i)}{\prod_{i=1}^{B  }(1 - q^i)
           \prod_{i=1}^{A-B}(1 - q^i)}
\end{equation}
for $0\le B\le A$, and $\left[ {A \atop B} \right]_q=0$ otherwise.


\subsection{Parafermions}\label{ParaPathSec}

In regime II of \cite{abf} (the parafermion model),
each path $h$ is assigned a weight $\bwt^{\rm II}(h)$ given by:
\begin{equation}\label{BwtIIDef}
\bwt^{\rm II}(h)=\sum_{i=1}^L ic^{\rm II}(h_{i-1},h_{i},h_{i+1}),
\end{equation}
where the function $c^{\rm II}(h_{i-1},h_{i},h_{i+1})$
is defined by:
\begin{eqnarray*}
c^{\rm II}(h-1,h,h+1)&=&0\,;\\
c^{\rm II}(h+1,h,h-1)&=&0\,;\\
c^{\rm II}(h-1,h,h-1)&=&1/2\,;\\
c^{\rm II}(h+1,h,h+1)&=&1/2\,.
\end{eqnarray*}

The generating function for paths in the parafermion model is defined to be:
\begin{equation}\label{GenIIDef}
x^{p'}_L(a,b,c;q)=\sum_{h\in{\P}^{p'}_{a,b,c}(L)} q^{\bwt^{\rm II}(h)}.
\end{equation}
Then define $x^{p'}_L(a,b,c)=x^{p'}_L(a,b,c;q)$. We immediately obtain:

\begin{lemma}\label{ParaABFLem}
\begin{equation}
x^{p'}_{L}(a,b,c;q)=q^{{1\over4}L(L+1)}X^{p'}_{L}(a,b,c;q^{-1}).
\end{equation}
\end{lemma}

\Proof For each $h\in{\P}^{p'}_{a,b,c}(L)$,
$$
\bwt^{\rm II}(h)+\bwt^{\rm III}(h)=\sum_{i=1}^L
i(c^{\rm II}(h_{i-1},h_{i},h_{i+1})+c^{\rm III}(h_{i-1},h_{i},h_{i+1})).
$$
Since for each $i$, one of $c^{\rm II}h_{i-1},h_{i},h_{i+1})$ and
$c^{\rm III}(h_{i-1},h_{i},h_{i+1})$ is $1/2$ and the other is 0,
we get
$$
\bwt^{\rm II}(h)+\bwt^{\rm III}(h)=\sum_{i=1}^L {i/2}=L(L+1)/4.
$$
The result now follows from (\ref{GenIIIDef}) and (\ref{GenIIDef}).
\cqfd
\medskip

\subsection{An alternative prescription for weights}\label{AltPresSec}

In this section, we define the weight of a path in yet another way.
However, as we will see, the difference between that given here
and that given above is just an overall factor.

The new definition of the weight involves the path picture.
Consider paths in the set ${\P}^{p'}_{a,b,c}(L)$,
and define new coordinates on the grid as follows:
$$
x={i-(h-a)\over2},\qquad y={i+(h-a)\over2}.
$$
Thus, the $xy$-coordinate system has its origin at the path's initial
point, and is slanted at $45^o$ to the original $ih$-coordinate system.
Note that at each step in the path, either $x$ or $y$ is incremented
and the other is constant. In this system, the path depicted
in Fig.~1 has its first few coordinates at
$(0,0)$, $(0,1)$, $(0,2)$, $(1,2)$, $(1,3)$, $(2,3)$, $(3,3)$, $(4,3)$,
$\ldots$.
Now, if the $i$th vertex has coordinates $(x,y)$,
we define $c(h_{i-1},h_i,h_{i+1})$
according to the shape of the vertex as follows:
\begin{eqnarray*}
c(h-1,h,h+1)&=&0\,;\\
c(h+1,h,h-1)&=&0\,;\\
c(h-1,h,h-1)&=&x\,;\\
c(h+1,h,h+1)&=&y\,.
\end{eqnarray*}
\noindent
We shall refer to those vertices for which, in general,
the contribution is non-zero, as {\em scoring} vertices.
The other vertices will be termed {\em non-scoring}.

We now define
\begin{equation}\label{WtDef}
\wt(h)=\sum_{i=1}^L c(h_{i-1},h_i,h_{i+1}).
\end{equation}

To illustrate this procedure, consider again the path $h$ depicted in
Fig.~1 and take $c=3$. The above table indicates that there are scoring
vertices at $i=2$, $3$, $4$, $7$, $8$, $10$.
This leads to
$$
\wt(h)=0+2+1+3+4+4=14.
$$

We now define the generating function:
\begin{equation}\label{GenDef}
\chi^{1,p'}_{a,b,c}(L;q)
=\sum_{h\in{\P}^{p'}_{a,b,c}(L)} q^{\wt(h)},
\end{equation}
and set $\chi^{1,p'}_{a,b,c}(L)=\chi^{1,p'}_{a,b,c}(L;q)$.

We note a symmetry of these generating functions:
\begin{lemma}\label{ReflectLem}
Let $L\ge0$, $1\le a,b<p'$ and $c=b\pm1$. Then
$$
\chi^{1,p'}_{a,b,c}(L)=\chi^{1,p'}_{p'-a,p'-b,p'-c}(L).
$$
\end{lemma}

\Proof Let $h'$ be the path obtained from $h$ by reflecting it
in a horizontal axis so that $h'_i=p'-h_i$.
We immediately see that $\wt(h')=\wt(h)$.
The lemma then follows from the definition (\ref{GenDef}).
\cqfd
\medskip

\begin{lemma}\label{RenormParaLem}
Let $L\ge0$ and $1\le a,b<p'$. Then
\begin{equation}
\chi^{1,p'}_{a,b,b\pm1}(L)
=q^{-{1\over4}(L\pm(a-b))} x^{p'}_{L}(a,b,b\pm1)
\end{equation}
\end{lemma}

\Proof
Let $h\in{\P}^{p'}_{a,b,c}(L)$ and let $h$ have $N$ scoring vertices.
Let the $i$ coordinates of these vertices be $i_1,i_2,\ldots,i_N$,
with $1\le i_N<i_{N-1}<\cdots<i_1\le L$.
Then let $(x_j,y_j)$ be the $(x,y)$-coordinates of the peak
at $(i_j,h_{i_j})$.

If $c=b+1$, then there is a peak-down vertex at $i=i_1$.
Then the $x$-coordinates of $(i_1,h_{i_1})$ and $(L,b)$ are
equal so that $x_1={1\over2}(L+(a-b))$.
Furthermore, $y_1=y_2$, $x_2=x_3$, etc., so that
$x_j=x_{j+1}$ for $j$ even, and $y_j=y_{j+1}$ for $j$ odd;
with finally $x_N=0$ if $N$ is even, and $y_N=0$ if $N$ is odd.
Thereupon:
\begin{eqnarray*}
i_1 +i_2 + \cdots + i_N \hskip-15mm&&\\
&=&
(x_1+y_1)+(x_2+y_2)+\cdots+(x_N+y_N)\\
&=&
{1 \over 2}(L+(a-b))+2y_1+2x_2+2y_3+
\cdots
+
\begin{cases}
2x_N & \text{if $N$ is even},\\
2y_N & \text{if $N$ is odd}
\end{cases} \\
&=&
{1\over2}(L+(a-b))+2\,\wt(h).
\end{eqnarray*}

On the other hand, if $c=b-1$ whereupon $y_1={1\over2}(L-(a-b))$,
a similar argument (in fact, by just exchanging the roles of $x$
and $y$ in the above) leads to:
$$
i_1+i_2+\cdots+i_N
\;=\;
{1\over2}(L-(a-b))+2\,\wt(h).
$$

Combining these two results thus yields:
$$
\wt^{\rm II}(h)={1\over4}(L\pm(a-b))+\wt(h),
$$
when $b=c\pm1$.
The lemma then follows immediately from the
definitions (\ref{GenIIDef}) and (\ref{GenDef}).
\cqfd
\medskip

The following result now provides the relationship between the
renormalised ABF and parafermion generating functions.

\begin{lemma}\label{DualityLem}
Let $L\ge0$, $1\le a,b<p'$ and $c=b\pm1$. Then
$$
\chi^{1,p'}_{a,b,c}(L;q)
=q^{{1\over4}(L^2-(a-b)^2)}\chi^{p'-1,p'}_{a,b,c}(L;q^{-1}).
$$
\end{lemma}

\Proof
This result follows from combining the definition (\ref{RenormABFEq})
with Lemmas \ref{ParaABFLem} and \ref{RenormParaLem}.
\cqfd
\medskip

\subsection{Striking sequence of a path}\label{StrikingSec}

Scanning from left to right, one can think of each $h\in{\P}^{p'}_{a,b}(L)$ 
as a sequence of straight lines, alternating in direction between NE and 
SE. Let the lengths of these lines be $w_1$, $w_2$, $w_3,\ldots,w_l,$ for 
some $l$, so that $w_1+w_2+\cdots+w_l=L(h)$. 
In what follows, we permit $w_1=0$ and $w_l=0$, but restrict
$w_i>0$ for $1<i<l$.

As will become clear shortly, we need to augment the definition of a path
as follows: for each path, we fix $e,f\in\{0,1\}$ (arbitrarily for now), 
and require $w_1$ to be the number of SE (resp.\ NE) segments at the 
beginning of the path if $e=0$ (resp.\ $e=1$), and $w_l$ to be the number 
of NE (resp.\ SE) segments at the end of the path if $f=0$ (resp.\ $f=1$). 
Notice that there are 4 possible augmentations of each path.
This definition implies that $l\equiv e+f\,(\mod2)$.
The striking sequence of the 
\lq augmented\rq\ path, $h$, is then defined to be the symbol:
\begin{equation}\label{HseqDef}
\left( w_1,w_2,w_3,\ldots,w_l \right)^{(e,f)}.
\end{equation}

We now define $m^{(e,f)}(h)=L-l+2$, where $l$ is the number of 
elements in the striking sequence above. Note that 
$l \equiv e+f\,(\mod2)$ implies that 
$m^{(e,f)}(h)\equiv L+e+f\,(\mod2)$.
For a given path $h$, we see that no two values of $m^{(0,0)}(h)$,
$m^{(0,1)}(h)$, $m^{(1,0)}(h)$ and $m^{(1,1)}(h)$ are guaranteed
equal.

For example, the path $h$ shown in Fig.~1
has the four possible striking sequences:
$(0,2,1,1,3,1,2,1)^{(0,0)}$, $(2,1,1,3,1,2,1)^{(1,0)}$,
$(0,2,1,1,3,1,2,1,0)^{(0,1)}$ and $(2,1,1,3,1,2,1,0)^{(1,1)}$.
Thence, we obtain $m^{(0,0)}(h)=5$, $m^{(1,0)}(h)=6$,
$m^{(0,1)}(h)=4$ and $m^{(1,1)}(h)=5$.


\begin{lemma}\label{CornerSwitchLem}
Let $h\in{\P}^{p'}_{a,b}(L)$ and let $e,f\in\{0,1\}$.
\begin{itemize}
\item $m^{(1,f)}(h)=m^{(0,f)}(h)+1$, for $a=1$;
\item $m^{(1,f)}(h)=m^{(0,f)}(h)-1$, for $a=p'-1$;
\item $m^{(e,1)}(h)=m^{(e,0)}(h)+1$, for $b=1$;
\item $m^{(e,1)}(h)=m^{(e,0)}(h)-1$, for $b=p'-1$.
\end{itemize}
\end{lemma}

\Proof If $a=1$, then the first segment of the path is certainly
in the NE direction.
Thus, with $w_1>0$, the path $h$ has striking sequences
$(w_1,w_2,\ldots,w_l)^{(1,f)}$ and $(0,w_1,w_2,\ldots,w_l)^{(0,f)}$.
Then $m^{(1,f)}=L-l+2$ and $m^{(0,f)}=L-(l+1)+2$,
whereupon the first result follows immediately.
The other three results follow in an analogous way.
\cqfd
\medskip

The purpose of assigning $e$ and $f$ to a path $h$ of length $L$,
is to enable the $0$th and $L$th vertices to be each designated as
scoring or non-scoring. In fact, if these vertices are included,
$m^{(e,f)}(h)$ gives the total number of non-scoring vertices in $h$.
We see that prescribing $e$ and $f$ is equivalent to appending 
two extra segments to the path, {\sl a pre-segment} that ends 
at the 0th vertex, and {\sl a post-segment} that starts at the 
Lth vertex. 
Setting $e=0$ (resp.\ $e=1$) is equivalent to having a  pre-segment that 
points SE (resp.\ NE). 
Setting $f=0$ (resp.\ $f=1$) is equivalent to having a post-segment that 
points NE (resp.\ SE). 

If we refer to these additional segments as the $0$th and $(L+1)$th
segments respectively, then the above definition of the striking
sequence implies that $w_1$ counts the number of segments
at the beginning of the path in the same direction as (but not including)
the $0$th and $w_l$ counts the number of segments at the end of
the path in the same direction as (but not including) the $(L+1)$th.
Note that specifying $f$ is equivalent to specifying the extra
point $h_{L+1}$.

We now define a weight for each $h\in{\P}^{p'}_{a,b}(L)$
that depends on $e$ and $f$,
and then show that this weight is equal to the weight of the
corresponding path in ${\P}^{p'}_{a,b,c}(L)$ with $c$
appropriately defined.

\begin{definition}\label{WtStrikeDef}
Let $h\in{\P}^{p'}_{a,b}(L)$, let $e,f\in\{0,1\}$ and
let the path $h$ have the striking sequence
$(w_1,w_2,w_3,\ldots,w_l)^{(e,f)}$.
Then define
$$
\wt^{(e,f)}(h)=\sum_{i=2}^{l-1} (w_{i-1}+w_{i-3}
+\cdots+w_{1+i\,\mbox{\scriptsize\rm mod}\,2}).
$$
\end{definition}

\noindent
We now show that this definition essentially provides the weight
of the corresponding path for which the appropriate extra
point $h_{L+1}$ is defined.

\begin{lemma}\label{WtStrikeLem}
Let $h\in{\P}^{p'}_{a,b}(L)$ and let $e,f\in\{0,1\}$.
If $f=0$ then let $c=b+1$ and if $f=1$ then let $c=b-1$.
Then let $h'\in{\P}^{p'}_{a,b,c}(L)$
be such that $h'_i=h_i$ for $0\le i\le L$.
Then $\wt(h')=\wt^{(e,f)}(h)$.
\end{lemma}

\Proof
Let the path $h$ have the striking sequence
$(w_1,w_2,w_3,\ldots,w_l)^{(e,f)}$,
where $w_1\ge0$, $w_l\ge0$ and $w_i>0$ for $1<i<l$.

Except for $i=l$ and possibly $i=1$, there is a scoring vertex at
the end of the $i$th line (which has length $w_i$) of $h'$.
First assume that the first $w_1$ segments of $h$ are in the NE
direction.
Then, for $i$ odd, the $i$th line is in the NE direction and its
$x$-coordinate is $w_2+w_4+\cdots+w_{i-1}$. By the prescription
of the previous section, this
line thus contributes $(w_2+w_4+\cdots+w_{i-1})$ to the weight
$\wt(h')$ of $h'$. Similarly, for $i$ even, the $i$th line
is in the SE direction and contributes
$(w_1+w_3+\cdots+w_{i-1})$ to $\wt(h')$.
This proves the lemma if the first $w_1$ segments are in the NE direction.
The reasoning is almost identical for the other case.
\cqfd
\medskip

\begin{lemma}\label{WtSwitchLem}
Let $h\in{\P}^{p'}_{a,b}(L)$ and $e,f\in\{0,1\}$.
Then $\wt^{(0,f)}(h)=\wt^{(1,f)}(h)$. Furthermore,
\begin{itemize}
\item $\wt^{(e,0)}(h)=\wt^{(e,1)}(h)+{1\over2}(L-a +1)$, for $b=1   $;
\item $\wt^{(e,0)}(h)=\wt^{(e,1)}(h)-{1\over2}(L-p'+a)$, for $b=p'-1$.
\end{itemize}
\end{lemma}

\Proof For either $e=0$ or $e=1$, $h$ has striking sequence
$(w_1,w_2,\ldots,w_l)^{(e,f)}$ with $w_1>0$.
Then $h$ also has striking sequence
$(0,w_1,w_2,\ldots,w_l)^{(1-e,f)}$.
Definition \ref{WtStrikeDef} then gives $\wt^{(e,f)}(h)=\wt^{(1-e,f)}(h)$,
thereby proving the first part.

If $b=1$, then the $L$th segment of $h$ is necessarily in the SE direction.
Therefore $h$ has striking sequence $(w_1,w_2,\ldots,w_l)^{(e,1)}$
with $w_l>0$. It then also has striking sequence
$(w_1,w_2,\ldots,w_l,0)^{(e,0)}$.
Therefore, $\wt^{(e,0)}(h)-\wt^{(e,1)}(h)=w_{l-1}+w_{l-3}+w_{l-5}+\cdots$.
Since $w_1+w_2+\cdots+w_l=L$ and
$(w_l+w_{l-2}+\cdots)-(w_{l-1}+w_{l-3}+\cdots)=a-1$,
it follows that $\wt^{(e,0)}(h)-\wt^{(e,1)}(h)={1\over2}(L-a+1)$.

A similar argument for the $b=p'-1$ case, shows that $h$ has striking
sequences $(w_1,w_2,\ldots,w_l)^{(e,0)}$
and $(w_1,w_2,\ldots,w_l,0)^{(e,1)}$, where $w_l>0$.
Then $\wt^{(e,1)}(h)-\wt^{(e,0)}(h)=w_{l-1}+w_{l-3}+w_{l-5}+\cdots$.
Now $(w_l+w_{l-2}+\cdots)-(w_{l-1}+w_{l-3}+\cdots)=p'-a+1$,
whereupon $\wt^{(e,1)}(h)-\wt^{(e,0)}(h)={1\over2}(L-p'+a)$.
\cqfd
\medskip

\subsection{Restricted generating functions}\label{RestGenSec}

We now define the set of paths ${\P}^{p'}_{a,b,e,f}(L,m)$
to be the subset of ${\P}^{p'}_{a,b}(L)$, comprising
those paths $h$ for which $m^{(e,f)}(h)=m$.
Let $\chi^{1,p'}_{a,b,e,f}(L,m)$ be the generating function for all
such paths:
$$
\chi^{1,p'}_{a,b,e,f}(L,m)=
\sum_{h\in{\P}^{p'}_{a,b,e,f}(L,m)} q^{\wt^{(e,f)}(h)}.
$$
Note that $\chi^{1,p'}_{a,b,e,f}(L,m)=0$ unless $m\equiv L+e+f\,(\mod2)$.

\begin{lemma}\label{GenSwitchLem}
Let $1\le a,b<p'$ and $e,f\in\{0,1\}$. Then:
\begin{itemize}
\item $\chi^{1,p'}_{1,b,1,f}(L,m)=\chi^{1,p'}_{1,b,0,f}(L,m-1)$;
\item $\chi^{1,p'}_{p'-1,b,1,f}(L,m)=\chi^{1,p'}_{p'-1,b,0,f}(L,m+1)$;
\item $\chi^{1,p'}_{a,1,e,0}(L,m)
 =q^{{1\over2}(L-a+1)}\chi^{1,p'}_{a,1,e,1}(L,m+1)$;
\item $\chi^{1,p'}_{a,p'-1,e,0}(L,m)
 =q^{-{1\over2}(L-p'+a)}\chi^{1,p'}_{a,p'-1,e,1}(L,m-1)$.
\end{itemize}
\end{lemma}

\Proof If $h\in{\P}^{p'}_{1,b,1,f}(L,m)$ then
$h\in{\P}^{p'}_{1,b}(L)$ and $m=m^{(1,f)}(h)$.
Then by Lemma \ref{CornerSwitchLem}, $m^{(0,f)}(h)=m-1$,
so that
${\P}^{p'}_{1,b,0,f}(L,m-1)\subset{\P}^{p'}_{1,b,1,f}(L,m)$.
On reversing the argument, the direction of the inclusion here
is changed, whereupon
${\P}^{p'}_{1,b,0,f}(L,m-1)={\P}^{p'}_{1,b,1,f}(L,m)$.
By Lemma \ref{WtSwitchLem}, $\wt^{(0,f)}(h)=\wt^{(1,f)}(h)$
for all $h\in{\P}^{p'}_{1,b}(L)$.
Thereupon, the first statement of the lemma is proved. The second
statement follows in a similar way.

For the third statement, we also obtain
${\P}^{p'}_{a,1,e,0}(L,m-1)={\P}^{p'}_{a,1,e,1}(L,m)$
in a similar way.
However, Lemma \ref{WtSwitchLem} gives
$\wt^{(e,0)}(h)=\wt^{(e,1)}(h)+{1\over2}(L-a+1)$
for all $h\in{\P}^{p'}_{a,1}(L)$.
Thereupon, the third statement follows.
The fourth statement follows similarly.
\cqfd
\medskip

The following result will act as a seed to generate further expressions.

\begin{lemma}\label{SeedLem}
Let $1\le a,b<p'$ with $p'\ge2$. In addition, let $e,f\in\{0,1\}$. Then
$$
\chi^{1,p'}_{a,b,e,f}(0,m)= \delta_{a,b}\delta_{\vert e-f\vert,m}.
$$
If $L\ge0$, then
$$
\chi^{1,2}_{1,1,e,f}(L,m)= \delta_{L,0}\delta_{\vert e-f\vert,m}.
$$

\end{lemma}

\Proof Clearly ${\P}^{p'}_{a,b}(0)=\emptyset$ if $a\ne b$.
Otherwise, it contains a single element.
Let $h$ designate this path.
In the case $e=f$, $h$ has striking sequence $(0,0)^{(e,f)}$
whereupon $m^{(e,f)}(h)=0$ immediately from the definition.
We also obtain $\wt^{(e,f)}(h)=1$ whereupon the result
follows for $e=f$. In the case $e\ne f$, $h'$ has striking
sequence $(0)^{(e,f)}$.
Then $m^{(e,f)}(h)=1$ and $\wt^{(e,f)}(h)=1$ whence the required
result also follows for $e\ne f$.

For the second expression, it is clear that
${\P}^{2}_{a,b}(L)=\emptyset$ unless $L=0$.
The result then follows from the first part.
\cqfd
\medskip

\begin{lemma}\label{FirstStepLem}
Let $1\le a,b<p'$ and $c=b\pm1$. Then if $c=b+1$ let $f=0$,
and if $c=b-1$ let $f=1$.
Then for each $e\in\{0,1\}$,
$$
\chi^{1,p'}_{a,b,c}(L)=\sum \chi^{1,p'}_{a,b,e,f}(L,m)
$$
where the sum is over all $m$ for which 
$m\equiv L+e+f\,(\mod2)$.
\end{lemma}

\Proof 
For each $h'\in{\P}^{p'}_{a,b,c}(L)$,
there is a corresponding path $h\in{\P}^{p'}_{a,b}(L)$
for which $h'_i=h_i$ for $0\le i\le L$, and vice-versa.
Moreover, $h\in{\P}^{p'}_{a,b,e,f}(L,m)$
for one and only one value of $m$ which is given by $m=m^{(e,f)}(h)$,
whereupon $m\equiv L+e+f\,(\mod2)$.
Then ${\P}^{p'}_{a,b,c}(L)=\bigcup_{m\equiv L+e+f\,
(\mbox{\scriptsize mod2})} {\P}^{p'}_{a,b,e,f}(L,m)$.
Since, by Lemma \ref{WtStrikeLem}, $\wt(h')=\wt^{(e,f)}(h)$,
the current lemma follows.
\cqfd
\medskip

\section{Path transformations}\label{TranSec}

\subsection{$\B$-transforms}\label{BTranSec}

In this section, we define a method of transforming a
path in ${\P}^{p'}_{a,b}(L)$ to yield one in
${\P}^{p'+1}_{a',b'}(L')$ for certain $a'$, $b'$ and $L'$.
This transform will be referred to as a
${\B}$-transform%
\footnote{It may be seen that when $e=f=0$, the ${\B}$-transform
described here is a generalisation of the $\B$-transform described in
\cite{flpw} as it acts upon paths in ${\P}^{1,p'}_{1,1}(L)$.}.

The action of the ${\B}$-transform is most easily described using
the striking sequences.
For $e,f\in\{0,1\}$, the action of the
${\B}$-transform on the path $h$ described by
the striking sequence $(w_1,w_2,w_3,\ldots,w_l)^{(e,f)}$
is to yield the path $\hat h$ with $\hat h_0=h_0+e$ and
striking sequence $(w_1,w_2+1,w_3+1,\ldots,w_{l-1}+1,w_l)^{(e,f)}$.
Note that the ${\B}$-transform action is dependent on the values
of $e$ and $f$ that appear in the striking sequence of the path.

For example, if $e=f=0$ then the action of the
${\B}$-transform on the path given in Fig.~1 results in the path:

\medskip
\centerline{\epsfig{file=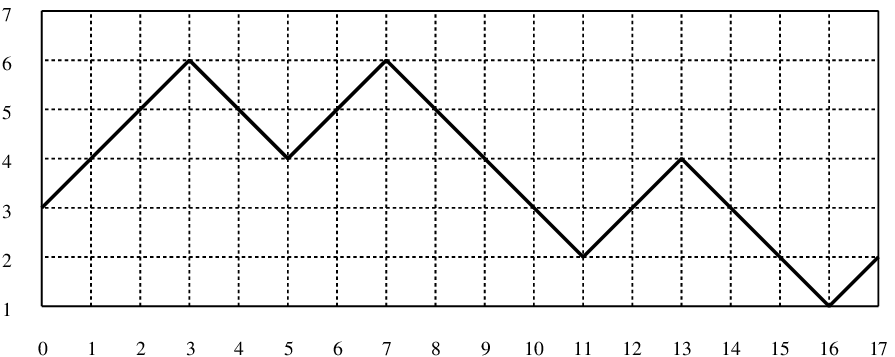}}
\medskip

\noindent
and if $e=1$ and $f=0$ then the action of the
${\B}$-transform on the path given in Fig.~1 results in the path:

\medskip
\centerline{\epsfig{file=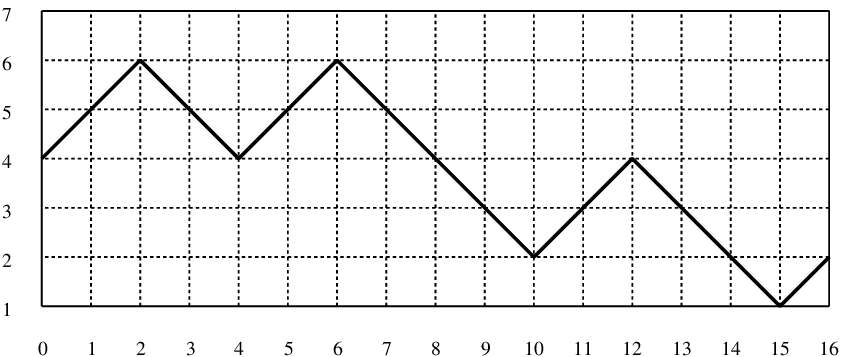}\;.}
\medskip

\noindent
Note that the path obtained from the action of the ${\B}$-transform
is such that there are no two consecutive scoring vertices.

Special care must be taken when dealing with paths of length 0
when $e\ne f$. The striking sequence of such a path is $(0)$.
We choose to leave the action of the $\B$-transform on such
a path undefined. Then Lemmas \ref{BHashLem}, \ref{BWtLem},
\ref{TranLem}, \ref{WtChngLem1}, \ref{WtChngLem2} and
\ref{GaussLem} will not apply for such paths.
However, they appear as a special case in the proof of
Lemma \ref{InductLem}.

\begin{lemma}\label{BHashLem}

Let $h\in{\P}^{p'}_{a,b}(L)$ and, for $e,f\in\{0,1\}$,
let $\hat h$ be the path obtained from
the action of the ${\B}$-transform on $h$.
Then $\hat h\in{\P}^{p'+1}_{a+e,b+f}(\hat L)$,
$m^{(e,f)}(\hat h)=L$ and
$L(\hat h)=\hat L=2L-m^{(e,f)}(h)$.
\end{lemma}

\Proof From the definition of the ${\B}$-transform, we immediately
obtain $\hat L=L+l-2$, whereupon $m^{(e,f)}(h)=L-l+2$
implies that $\hat L=2L-m^{(e,f)}(h)$.
Additionally $m^{(e,f)}(\hat h)=\hat L-l+2=L$.
Now set $\hat a=\hat h_0$ and $\hat b=\hat h_{\hat L}$.
The definition of the ${\B}$-transform immediately gives
$\hat a=a+e$.
In terms of the striking sequence $(w_1,w_2,\ldots,w_l)^{(e,f)}$
of $h$, we have $a-b=(-1)^e((w_1+w_3+\cdots)-(w_2+w_4+\cdots))$.
Then if $l$ is even, $\hat a-\hat b=a-b$, and if $l$ is odd
$\hat a-\hat b=a-b-(-1)^e$. 
Using $l\equiv e+f\,(\mod2)$ gives $\hat a-\hat b=a-b+e-f$
in both cases and hence $\hat b=b+f$ as required.
\cqfd
\medskip

\begin{lemma}\label{BWtLem}
Let $h\in{\P}^{p'}_{a,b}(L)$ and, for $e,f\in\{0,1\}$,
let $\hat h$ be the path obtained from
the action of the ${\B}$-transform on $h$.
Then,
$$
\wt^{(e,f)}(\hat h)=
\wt^{(e,f)}(h)+{1\over4}((\hat L-\hat m)^2-\delta_{e+f,1}),
$$
where $\hat L=L(\hat h)$ and $\hat m=m^{(e,f)}(\hat h)$.
\end{lemma}

\Proof Let the striking sequence of $h$ be
$(w_1,w_2,\ldots,w_{l-1},w_l)^{(e,f)}$ whereupon that of $\hat h$ is
$(w_1,w_2+1,\ldots,w_{l-1}+1,w_l)^{(e,f)}$.
Definition \ref{WtStrikeLem} then gives:
$$
\wt^{(e,f)}(\hat h)-\wt^{(e,f)}(h)=
0+1+1+2+2+3+\cdots+
\lfloor(l-3)/2\rfloor+
\lfloor(l-2)/2\rfloor.
$$
This sum is ${1\over4}(l-2)^2$ if $l$ is even
and ${1\over4}((l-2)^2-1)$ if $l$ is odd.
The result then follows because $L(\hat h)-m^{(e,f)}(\hat h)=l-2$
and $l\equiv e+f\,(\mod2)$.
\cqfd
\medskip

\subsection{Inserting particles}\label{InsertSec}

Given a path $h^{(0)}$ of length $L$ and striking sequence
$(w_1,w_2,\ldots,w_l)^{(e,f)}$, we may extend $h^{(0)}$ by a
process we refer to as {\it inserting particles}.
If $a=1$, we restrict this process to the $e=0$ case,
and if $a=p'-1$, we restrict to the $e=1$ case.
The effect of inserting one particle is to produce
a path $h^{(1)}$ with the same starting point and striking sequence
$(0,1,w_1+1,w_2,\ldots,w_l)^{(e,f)}$. Thus the new path
has length $L+2$.
Notice that the way that the path is extended depends on $e$.
By iterating the process, we may insert $k$ particles
into $h^{(0)}$ to obtain a path $h^{(k)}$ of length $L+2k$.

\begin{lemma}\label{TranLem}
Let $h\in{\P}^{p'}_{a,b}(L)$ and, for $e,f\in\{0,1\}$,
let $h^{(0)}$ be the path obtained from
the action of the ${\B}$-transform on $h$,
and obtain $h^{(k)}$ from $h^{(0)}$ by inserting $k$ particles.
If $m=m^{(e,f)}(h)$, $m^\prime=m^{(e,f)}(h^{(k)})$ and
$L^\prime=L(h^{(k)})$, then $m'=L$,
$$
L'+m=2m'+2k,
$$
and
$$
\wt^{(e,f)}(h^{(k)})=
\wt^{(e,f)}(h)+{1\over4}\left((L'-m')^2-\delta_{e+f,1}\right).
$$
\end{lemma}

\Proof
By Lemma \ref{BHashLem}, $L(h^{(0)})=2L-m$ and $m^{(e,f)}(h^{(0)})=L$.
Inserting $k$ particles then gives $L'=L(h^{(0)})+2k=2L-m+2k$
and $m'=m^{(e,f)}(h^{(k)})=m^{(e,f)}(h^{(0)})=L$.

To obtain the final result, let $h^{(0)}$ have striking sequence
$(w_1,w_2,\ldots,w_l)^{(e,f)}$,
whereupon that of $h^{(1)}$ is
$(0,1,w_1+1,w_2,\ldots,w_l)^{(e,f)}$.
Then, $m^{(e,f)}(h^{(1)})=m^{(e,f)}(h^{(0)})$
and Definition \ref{WtStrikeLem} gives
$\wt^{(e,f)}(h^{(1)})=\wt^{(e,f)}(h^{(0)})+l-1$.
Repeated application then yields
$m^{(e,f)}(h^{(k)})=m^{(e,f)}(h^{(0)})$ and
\begin{eqnarray*}
\wt^{(e,f)}(h^{(k)})&=&
\wt^{(e,f)}(h^{(0)})+k(l-1)+k(k-1)\\
&=&\wt^{(e,f)}(h^{(0)})+k(l-2)+k^2\\
&=&\wt^{(e,f)}(h^{(0)})+k\left(L(h^{(0)})-m^{(e,f)}(h^{(0)})\right)+k^2\\
&=&\wt^{(e,f)}(h)+{1\over4}\left(\left(L(h^{(0)})
-m^{(e,f)}(h^{(0)})\right)^2-\delta_{e+f,1}\right)\\
&&\qquad\qquad\qquad\qquad
{}+k\left(L(h^{(0)})-m^{(e,f)}(h^{(0)})\right)+k^2,
\end{eqnarray*}
where the final equality follows from Lemma \ref{BWtLem}.
The required expression now results because, from above,
$L(h^{(0)})=L'-2k$ and $m^{(e,f)}(h^{(0)})=m'$.
\cqfd
\medskip

\subsection{Moving particles}\label{MovesSec}

In this section, we specify two types of local deformation of a path.
These deformations will be known as {\em moves}.
In each case, a particular sequence of four segments
of a path is changed to a different sequence, the remainder of the
path being unchanged.
The moves are as follows --- the path portion to the left of the arrow
is changed to that on the right:

\bigskip
\centerline{\epsfig{file=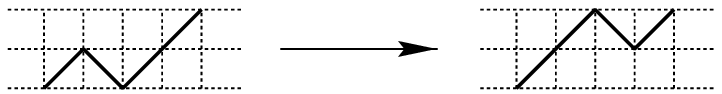}}
\nobreak
\centerline{Move.\ 1.}
\bigskip
\centerline{\epsfig{file=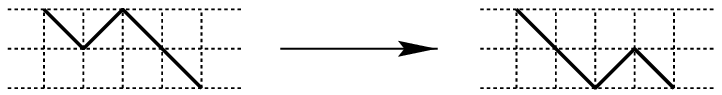}}
\nobreak
\centerline{Move.\ 2.}
\bigskip

\noindent Note that the two moves are inversions of one another.

\begin{lemma}\label{WtChngLem1}
Let $h$ be a path for which four consecutive segments
are as in one of diagrams on the left above.
Let $\hat h$ be that obtained from $h$ by changing those segments
according to the move. Then, for $e,f\in\{0,1\}$,
$$
\wt^{(e,f)}(\hat h)=\wt^{(e,f)}(h)+1.
$$
Additionally, $m^{(e,f)}(\hat h)=m^{(e,f)}(h)$
and $L(\hat h)=L(h)$.
\end{lemma}
\Proof For each case, take the $xy$-coordinate of the
leftmost point of this portion of a path to be $(x_0,y_0)$.
Now consider the contribution to the weight of the three
vertices in question before and after the move
In both cases, the contribution is
$x_0+y_0+1$ before the move and $x_0+y_0+2$ afterwards.
Thus the first result holds.
The final result follows immediately from the definitions.
\cqfd
\medskip

Now observe that for each of the moves specified above,
the sequence of path segments before the move
consists of an adjacent pair of scoring vertices followed by a
non-scoring vertex.
The specified move then consists of replacing such a combination
with a non-scoring vertex followed by two scoring vertices.
It is useful to interpret this as
the pair of adjacent scoring vertices having moved by one step.
In fact, it is useful to refer to a pair of adjacent scoring
vertices as a particle.
Thus both of the moves consists of a particle moving rightwards
by one step.

In addition to the moves described above, and depending on
the values of $e$ and $f$, we permit certain deformations
of a path close to its left and right extremities.
They are as follows.

\bigskip
If $e=1:\qquad\qquad\qquad
\raisebox{-10pt}[0pt]{\epsfig{file=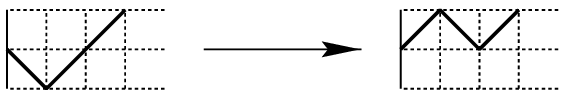}}$
\bigskip

\bigskip
If $e=0:\qquad\qquad\qquad
\raisebox{-10pt}[0pt]{\epsfig{file=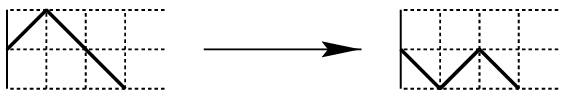}}$
\bigskip

\bigskip
If $f=0:\qquad\qquad\qquad
\raisebox{-10pt}[0pt]{\epsfig{file=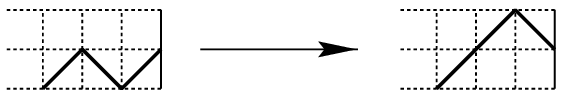}}$
\bigskip

\bigskip
If $f=1:\qquad\qquad\qquad
\raisebox{-10pt}[0pt]{\epsfig{file=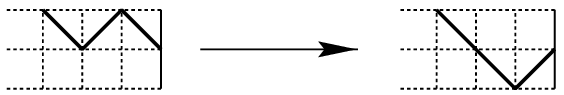}}$
\bigskip

In fact, the above four moves may be considered as instances of the
two moves described beforehand, if we append two extra segments
to the path as described in the paragraph
following Lemma \ref{CornerSwitchLem}.

\begin{lemma}\label{WtChngLem2}
Let $h$ be a path for which three consecutive segments
occupy three extreme positions and are as in one of
the diagrams above, on the left.
Let $\hat h$ be that obtained from $h$ by changing those segments
according to the move. Then, for $e,f\in\{0,1\}$,
$$
\wt^{(e,f)}(\hat h)=\wt^{(e,f)}(h)+1.
$$
Additionally, $m^{(e,f)}(\hat h)=m^{(e,f)}(h)$
and $L(\hat h)=L(h)$.
\end{lemma}
\Proof As for Lemma \ref{WtChngLem1}.
\cqfd
\medskip

With the interpretation of the extremity moves in terms of additional
segments to the left and right of the original path,
we may also designate the 0th and $L$th vertices as scoring or
non-scoring. We then see that the four extremity moves also
involve the replacing of a pair of scoring vertices that are
followed by a non-scoring vertex, by a non-scoring vertex
followed by a pair of scoring vertices.

\subsection{Waves of particles}\label{WavesSec}

Since each of the moves described above involves a pair of scoring
vertices moving rightwards by one step, we see that a succession of
such moves is possible until the pair is followed by another
scoring vertex.
If this itself is followed by yet another scoring vertex,
we forbid further movement. However, if it is followed by
a non-scoring vertex, further movement is allowed after
considering the latter two of the three consecutive scoring
vertices to be the particle (instead of the first two).

Now consider a path $h^{(k)}$ that results from inserting
$k$ particles into $h^{(0)}$ that itself results from
the action of the $\B$-transform.
We now show that the paths that result from moving these particles
in all possible ways, are indexed by partitions
(a definition of a partition may be found in Appendix B).

\begin{lemma}\label{GaussLem} 
There is a bijection between the paths obtained by moving the particles in
$h^{(k)}$ and the partitions $\lambda$ with at most $k$
parts that have no part larger than $m=m^{(e,f)}(h^{(k)})$.
This bijection is such that if
$h$ is the bijective image of a particular $\lambda$ then
$$
\wt(h)=\wt(h^{(k)})+\wt(\lambda),
$$
where $\wt(\lambda)=\lambda_1+\lambda_2+\cdots+\lambda_k$.
\end{lemma}

\Proof
If $h^{(0)}$ has striking sequence
$(w_1,w_2,\ldots,w_l)^{(e,f)}$ then there are $l-1$ scoring vertices
between the 0th and $L$th vertex inclusive --- one at the end
of the $i$th line which has length $w_i$ for $1\le i<l$.
Since there are altogether $L+1$ vertices, the number
of non-scoring vertices is $L+1-(l-1)=L-l+2=m^{(e,f)}(h^{(0)})$.

Since each particle moves by traversing a non-scoring vertex,
and there are $m$ of these to the right of the rightmost particle 
in $h^{(k)}$, and there are no consecutive scoring vertices to its 
right, this particle can make $\lambda_1$ moves to the right, with 
$0\le\lambda_1\le m$. Similarly, the next rightmost particle can 
make $\lambda_2$ moves to the right with $0\le\lambda_2\le\lambda_1$.
Here, the upper restriction arises because the two scoring vertices 
would then be adjacent to those of the first particle. Continuing in 
this way, we obtain that all possible final positions of the particles 
are indexed by $\lambda=(\lambda_1,\lambda_2,\ldots,\lambda_k)$ with 
$m\ge\lambda_1\ge\lambda_2\ge\cdots\ge\lambda_k\ge0$, that is, by 
partitions of at most $k$ parts with no part exceeding $m$. Moreover, 
since by Lemmas \ref{WtChngLem1} and \ref{WtChngLem2}
the weight increases by one for each 
move, the weight increase after the sequence of moves specified by 
a particular $\lambda$ is equal to $\wt(\lambda)$.
\cqfd
\medskip

We now combine the above results on the $\B$-transform,
and the inserting and moving of particles to obtain
an expression which enables the restricted generating
function $\chi^{1,p'}_{a,b,e,f}(L',m')$ to be obtained
from those of a \lq simpler\rq\ model.

\begin{lemma}\label{InductLem}
Let $e,f\in\{0,1\}$. Then, if $p'>2$, and $1\le a-e,b-f<p'-1$,
\begin{displaymath}
\chi^{1,p'}_{a,b,e,f}(L',m')
=
\hskip-3mm
\sum_{m\equiv L'\,(\mbox{\scriptsize\rm mod}\,2)}
\hskip-3mm
q^{{1\over4}\left((L'-m')^2-\delta_{e+f,1}\right)}
\left[{{1\over2}(L'+m)\atop m'}\right]_q
\chi^{1,p'-1}_{a-e,b-f,e,f}(m',m).
\end{displaymath}
\end{lemma}

\Proof For the moment, exclude the case where $m'=0$ and $L'$ is odd.
Consider a path $h$ that contributes to
$\chi^{1,p'-1}_{a-e,b-f,e,f}(m',m)$ so that
$L(h)=m'$ and $m^{(e,f)}(h)=m$.
Let $h^{(0)}$ result from the action of the
${\B}$-transform on $h$.
By Lemma \ref{TranLem}, inserting $k$ particles into
$h^{(0)}$ results in a path $h^{(k)}$ of length $L'$ if and only if
$2k=L'-2m'+m$.
(In particular, $m\equiv L'\,(\mod2)$.)
The generating function for all paths obtained by moving the
particles in $h^{(k)}$ is then given by
$$
q^{{1\over4}\left((L'-m')^2-\delta_{e+f,1}\right)}
\left[{k+m'\atop m'}\right]_q
q^{\wt(h)},
$$
on using the last expression of Lemma \ref{TranLem}, and using
Lemma \ref{GaussLem} after noting that
the Gaussian term $\left[{k+m'\atop m'}\right]$ is the generating
function for all partitions $\lambda$ with $k$ non-negative parts,
none exceeding $m'$.
Summing this expression over all
$h\in{\P}^{p'-1}_{a-e,b-f,e,f}(m',m)$ and
then all $m\equiv L'\,(\mod2)$,
results in the expression on the right side of that in the premise.

We now show that every path $h'\in{\P}^{p'}_{a,b,e,f}(L',m')$
arises through the above procedure and moreover, arises in a unique way.
So let $h'\in{\P}^{p'}_{a,b,e,f}(L',m')$.
Locate the leftmost pair of consecutive scoring vertices in $h'$,
and move them leftward by reversing the particle moves,
until they occupy the 0th and 1st positions.
Now ignoring these two vertices, do the same with the
next leftmost pair of consecutive scoring vertices,
moving them leftward until they occupy the 2nd and 3rd
positions.
Continue in this way until all consecutive scoring vertices
occupy the leftmost positions of the path.
Say there are $\kappa$ of them, and let $k=\lfloor \kappa/2\rfloor$.
Removing the first $2k$ segments then produces a path
$h^{(0)}$ which, by construction, has no pair of
consecutive scoring vertices.
Thus, its striking sequence will be of the form
$(w_1,w_2,\ldots,w_{l-1},w_l)^{(e,f)}$, with $w_i\ge2$ for
$1<i<l$.
Thus $h^{(0)}$ arises from the action of the
${\B}$-transform on the path $h$ that has striking
sequence $(w_1,w_2-1,\ldots,w_{l-1}-1,w_l)^{(e,f)}$,
and $h_0=h'_0-e$.
Since $h$ and $k$ are thereby determined uniquely,
the lemma is proved except when $m'=0$ and $L'$ odd%
\footnote{The above proof fails in this case because,
the only non-zero term on the right side has $m'=0$ and
$m=1$, and the ${\B}$-transform is not defined on such paths.
Moreover, in the second paragraph, we would obtain
$\kappa=L'+1$ and removing this number of segments from a path
of length $L'$ is clearly nonsensical.}.

In this exceptional case, there is only one path $h'$ in
${\P}^{p'}_{a,b,e,f}(L',m')$ since each of its $L'+1$ vertices
are scoring, and thus $h'$ has scoring sequence
$(1,1,1,\ldots,1)^{(e,f)}$. Lemma \ref{WtStrikeLem}
gives, via Definition \ref{WtStrikeDef},
$\wt(h')={1\over4}({L'}^2-1)$.
The form of the striking sequence also ensures that
$a-b=\pm1$, and moreover, $a-e=b-f$, so that $\vert e-f\vert=1$.
In the sum on the right side, for other than $m=1$, we have
$\chi^{1,p'-1}_{a-e,b-f,e,f}(0,m)=0$.
However, by Lemma \ref{SeedLem},
$\chi^{1,p'-1}_{a-e,b-f,e,f}(0,1)=1$.
The proof is then complete.
\cqfd
\medskip

\section{Parafermion generating functions}\label{ParaFunSec}

In the previous section, it was indicated that all paths in
${\P}^{p'}_{a,b}(L)$ can be obtained by applying a sequence 
of ${\B}$-transforms to ${\P}^{2}_{1,1}(L)$. We do this by 
repeated application of Lemma \ref{InductLem}. In fact, this 
approach leads to four different concise constant-sign
expressions for the path generating functions. We derive two 
of these in Sections \ref{FirstSystemSec} and \ref{SecondSystemSec}.

These two expressions are equal. Combinatorially, this follows 
from the fact that both count the same set of objects. They 
look different, because we interpret these objects differently.
An analytic proof of their equality can be found in \cite{schilling}.

The two further constant-sign expressions
result from applying the symmetry obtained in Lemma \ref{ReflectLem}.
Although we do not show this, these expressions may also be obtained
by using the ${\B}$-transforms in a way similar to which the first
two were derived.

The expressions that we obtain involve the Cartan matrix
$\boldC=\boldC^{(t)}$ of the finite
dimensional Lie algebra $A_t$, i.e.\ $\boldC^{(t)}$ is the
$t\times t$ tri-diagonal matrix with entries $\boldC_{ij}$ for
$0\le i,j\le t-1$ where, when the indices are in this range,
\begin{equation}\label{CartanDef}
\begin{array}{cccl}
\boldC_{j,j-1}=-1,
&\boldC_{j,j}=2,
&\boldC_{j,j+1}=-1,
&\hbox{for $j=0,1,\ldots,t-1$.}
\end{array}
\end{equation}

The expressions also involve $(p'-2)$-dimensional vectors ${\boldu}_{a,b}$
which depend on $a$ and $b$.
Define
${\boldu}_{a,b}=(u_1,u_2,\ldots,u_{p'-2})$
%
%
%
%
%
as follows:
\begin{equation}
\begin{array}{ll}
u_i=\delta_{i,a-1}+\delta_{i,b-1},
\end{array}
\end{equation}
for $1\le i\le p'-2$.

Also, given a $t$-dimensional vector ${\boldu}=(u_1,u_2,\ldots,u_t)$,
we define $(t-1)$-dimensional vectors
$\boldQ(\boldu)=(Q_1,Q_2,\ldots,Q_{t-1})$ and
$\boldR(\boldu)=(R_1,R_2,\ldots,R_{t-1})$ as follows%
\footnote{What we denote here as $\boldQ(\boldu_{a,b})$, is denoted
$\tilde{\boldQ}_{p'-a,p'-b}$ in (3.3) of \cite{melzer}, and as
${\boldQ}_{p'-a,p'-b}$ in (2.2) of \cite{warnaar2}.
What we denote here as $\boldR(\boldu_{a,b})$, is denoted
${\boldQ}_{a+1,p'-b}$ in (3.3) of \cite{melzer}, and as
${\boldR}_{a,p'-b}$ in (2.5) of \cite{warnaar2}.
}:
\begin{equation}
\begin{array}{ll}
Q_i=(u_{i+1}+u_{i+3}+u_{i+5}+\cdots)\;\mod2&1\le i< t;\\
R_i=(t-i+u_{i+1}+u_{i+3}+u_{i+5}+\cdots)\;\mod2&1\le i< t,
\end{array}
\end{equation}
where $u_i=0$ for $i>t$.\footnote{Note that
$\boldR(\boldu)=\boldQ(\boldu+\bolde_t)$ where the
$t$-dimensional $\bolde_t=(0,0,\ldots,0,1)$.}
In the expressions, obtained below,
we require summations over vectors $\boldm=(m_1,m_2,\ldots,m_{t-1})$
for which $m_i\equiv Q_i\,(\mod2)$ for $1\le i< t$.
We shall denote such a restriction on $\boldm$ by simply
$\boldm\equiv\boldQ(\boldu)$.

To illustrate these definitions, let
$p'=14$, $a=4$ and $b=8$. Then
$\boldu_{4,8}=(0,0,1,0,0,0,1,0,0,0,0,0)$,
$\boldQ(\boldu_{4,8})=(0,0,0,1,0,1,0,0,0,0,0)$ and
$\boldR(\boldu_{4,8})=(1,0,1,1,1,1,1,0,1,0,1)$.

\subsection{First system}\label{FirstSystemSec}

In this section, we consider a sequence of ${\B}$-transforms
governed by the following values:
\begin{equation}\label{EFEq1}
\begin{array}{lcl}
e_i=1&\rm for&1\le i<a;\\
e_i=0&\rm for&a\le i\le p'-1;\\
f_i=1&\rm for&1\le i<b;\\
f_i=0&\rm for&b\le i\le p'-1.
\end{array}
\end{equation}
It will be useful to define:

\begin{eqnarray*}
a_j &=& a-\sum_{i=1}^{j-1} e_i=
\begin{cases}
1     & \text{if $a\le j\le p'-1$}, \\
a-j+1 & \text{if $1\le j\le a   $}
\end{cases} \\
b_j &=& b-\sum_{i=1}^{j-1} f_i=
\begin{cases}
1     & \text{if $b \le j \le p'-1$}, \\ 
b-j+1 & \text{if $1 \le j \le b   $}
\end{cases}
\end{eqnarray*}
In fact, although we make no use of this, these values actually
give the start and endpoints of the paths in
the sequence of sets of paths that we generate.

\begin{theorem}\label{FermThrm1}
Let $1\le a,b<p'$ and $c=b-1$.
Then, with $\boldC=\boldC^{(p'-2)}$, $m_0=L$, and $m_{p'-2}=0$,
\begin{displaymath}
\begin{array}{l}
\displaystyle
\chi^{1,p'}_{a,b,c}(L)\\
\hskip5mm
\displaystyle
{}={}
q^{{1\over4}(a-b-L^2)}
\sum
q^{{1\over4}\tilde{\sboldm}\sboldC\tilde{\sboldm}^T-{1\over2}m_{b-1}}
\prod_{i=1}^{p'-3}
\left[{{1\over2}(m_{i-1}+m_{i+1}+\delta_{i,a-1}+\delta_{i,b-1})
\atop m_i}\right]_q\!,
\end{array}
\end{displaymath}
where the above sum is over all $\boldm=(m_1,m_2,\ldots,m_{p'-3})$
for which $\boldm\equiv\boldQ(\boldu_{a,b})$,
with $\tilde{\boldm}=(m_0,m_1,m_2,\ldots,m_{p'-3})$.
\end{theorem}

\Proof Let $t=p'-2$. Lemma \ref{InductLem} implies that:\footnote{In
this and all subsequent proofs, we take the symbol \lq$\equiv$\rq\ to
mean equivalence modulo 2.}
\begin{displaymath}
\begin{array}{l}
\displaystyle
\chi^{1,p'+1-i}_{a_i,b_i,e_i,f_i}(m_{i-1},m_i)\\
\hskip8mm
\displaystyle
{}={}
\sum_{m\equiv m_{i-1}}
q^{{1\over4}\left((m_{i-1}-m_i)^2-\delta_{e_i+f_i,1}\right)}
\left[{{1\over2}(m_{i-1}+m)\atop m_i}\right]_q
\chi^{1,p'-i}_{a_{i+1},b_{i+1},e_i,f_i}(m_{i},m)
\end{array}
\end{displaymath}
for $i=1,2,\ldots,t$. In the $i$th case, we replace the summation
variable $m$ with $m=m_{i+1}+\delta_{i,a-1}+\delta_{i,b-1}$.
Thereupon, $m_{i-2}\equiv m_i+\delta_{i,a}+\delta_{i,b}$
for $2\le i\le t+1$.
We now express $\chi^{1,p'-i}_{a_{i+1},b_{i+1},e_i,f_i}(m_{i},m)$
in terms of
$\chi^{1,p'-i}_{a_{i+1},b_{i+1},e_{i+1},f_{i+1}}(m_{i},m_{i+1})$.
When $i\ne a-1$ and $i\ne b-1$, we immediately obtain
$$
\chi^{1,p'-i}_{a_{i+1},b_{i+1},e_i,f_i}(m_{i},m_{i+1})
=\chi^{1,p'-i}_{a_{i+1},b_{i+1},e_{i+1},f_{i+1}}(m_{i},m_{i+1}),
$$
from (\ref{EFEq1}).
When $i=a-1$, so that $a_{i+1}=1$, $e_i=1$ and $e_{i+1}=0$,
Lemma \ref{GenSwitchLem} yields:
$$
\chi^{1,p'-i}_{a_{i+1},b_{i+1},e_i,f_i}(m_{i},m_{i+1}+1+\delta_{a,b})
=\chi^{1,p'-i}_{a_{i+1},b_{i+1},e_{i+1},f_{i+1}+\delta_{a,b}}
                                       (m_{i},m_{i+1}+\delta_{a,b}),
$$
and when $i=b-1$, so that $b_{i+1}=1$, $f_i=1$ and $f_{i+1}=0$,
Lemma \ref{GenSwitchLem} yields:
$$
\chi^{1,p'-i}_{a_{i+1},b_{i+1},e_i-\delta_{a,b},f_i}(m_{i},m_{i+1}+1)
=q^{{1\over2}(a_b-1-m_{b-1})}
\chi^{1,p'-i}_{a_{i+1},b_{i+1},e_{i+1},f_{i+1}}(m_{i},m_{i+1}).
$$
By Lemma \ref{SeedLem},
$\chi^{1,2}_{1,1,0,0}(m_{t},m_{t+1})=\delta_{m_t,0}\delta_{m_{t+1},0}$,
so that we require $m_{t}=m_{t+1}=0$ for a non-zero contribution,
whereupon the Gaussian polynomial in the $t$th summation is equal to 1.
Moreover, it then follows from 
$m_{i-2}\equiv m_i+\delta_{i,a}+\delta_{i,b}$ that
$(m_1,m_2,m_3,\ldots,m_{t-1})\equiv\boldQ(\boldu_{a,b})$.
We now calculate
${1\over2}(a_b-1)-{1\over4}\sum_{i=1}^{t}\delta_{e_i+f_i,1}={1\over4}(a-b)$,
by using $e_i+f_i=1$ in $\vert a-b\vert$ cases and
$a_b=1$ if $b\ge a$ and $a_b=a-b+1$ if $b\le a$.
Combining all the above yields:
\begin{displaymath}
\begin{array}{l}
\displaystyle
\chi^{1,p'}_{a_1,b_1,e_1,f_1}(m_0,m_1)\\
\hskip2mm
\displaystyle
{}={}
\hskip-1mm
\sum
q^{{1\over4}\left(\sum_{i=1}^t(m_{i-1}-m_i)^2+a-b\right)-{1\over2}m_{b-1}}
\prod_{i=1}^{t-1}
\left[\!{{1\over2}(m_{i-1}\!+\!m_{i+1}\!+\!\delta_{i,a-1}\!+\!\delta_{i,b-1})
\atop m_i}\!\right]_q\!,
\end{array}
\end{displaymath}
where $m_t=0$ and the sum is over all
$(m_2,m_3,\ldots,m_{t-1})\equiv(Q_2,Q_3,\ldots,Q_{t-1})$,
when $\boldQ(\boldu_{a,b})=(Q_1,Q_2,\ldots,Q_{t-1})$.
{}From $m_1\equiv Q_1$, we obtain $m_0+e_1+f_1\equiv Q_1$.
The theorem now follows after noting that
$$
\sum_{i=1}^t(m_{i-1}-m_i)^2=\tilde{\boldm}\boldC^{(t)}\tilde{\boldm}^T-L^2,
$$
and using Lemma \ref{FirstStepLem} in the form:
$$
\chi^{1,p'}_{a,b,c}(m_0)=\sum_{m_1\equiv Q_1}
\chi^{1,p'}_{a_1,b_1,e_1,f_1}(m_0,m_1),
$$
and noting that $a=a_1$ and $b=b_1$.
\cqfd
\medskip

A further expression is obtained by using the reflection symmetry
identified in Lemma \ref{ReflectLem}:

\begin{corollary}\label{FermCorol1}
Let $1\le a,b<p'$ and $c=b+1$.
Then, with $\boldC=\boldC^{(p'-2)}$, $m_0=L$, and $m_{p'-2}=0$,
\begin{displaymath}
\begin{array}{l}
\displaystyle
\chi^{1,p'}_{a,b,c}(L)
=q^{{1\over4}(b-a-L^2)} \\
\displaystyle
\hskip8mm
\times\sum
q^{{1\over4}\tilde{\sboldm}\sboldC\tilde{\sboldm}^T-{1\over2}m_{p'-b-1}}
\prod_{i=1}^{p'-3}
\left[{{1\over2}(m_{i-1}+m_{i+1}+\delta_{i,p'-a-1}+\delta_{i,p'-b-1})
\atop m_i}\right]_q\!,
\end{array}
\end{displaymath}
where the sum is over all $\boldm=(m_1,m_2,\ldots,m_{p'-3})$
for which $\boldm\equiv\boldQ(\boldu_{p'-a,p'-b})$,
with $\tilde{\boldm}=(m_0,m_1,m_2,\ldots,m_{p'-3})$.
\end{corollary}

\subsection{Second system}\label{SecondSystemSec}

In this section, we consider a sequence of ${\B}$-transforms
different to that used in the previous section.
This leads to a constant-sign expression for $\chi^{1,p'}_{a,b,b-1}(L)$
that differs from that obtained in Theorem \ref{FermThrm1}.

In this section we use:

\begin{equation}\label{EFEq2}
\begin{array}{lcl}
e_i=0&\rm for&1\le i< p'-a;\\
e_i=1&\rm for&p'-a\le i\le p'-1;\\
f_i=1&\rm for&1\le i<b;\\
f_i=0&\rm for&b\le i\le p'-1.
\end{array}
\end{equation}
We define:
\begin{eqnarray*}
a_j  &=&  a-\sum_{i=1}^{j-1} e_i=
\begin{cases}
p'-j & \text{if $p'-a\le j\le p'-1$}, \\
a    & \text{if $1\le j\le p'-a$}     
\end{cases}    \\
b_j  &=&  b-\sum_{i=1}^{j-1} f_i=
\begin{cases}
1    & \text{if $b\le j\le p'-1$}, \\
b-j+1& \text{if $1\le j\le b$}
\end{cases}
\end{eqnarray*}

\begin{theorem}\label{FermThrm2}
Let $1\le a,b<p'$ and $c=b-1$.
Then, with $\boldC=\boldC^{(p'-2)}$, $m_0=L$, and $m_{p'-2}=0$,
\begin{displaymath}
\begin{array}{l}
\displaystyle
\chi^{1,p'}_{a,b,c}(L)\\
\hskip2mm
\displaystyle
{}={}
q^{{1\over4}(a-b-L^2)}
\hskip-1mm
\sum
q^{{1\over4}\tilde{\sboldm}\sboldC\tilde{\sboldm}^T-{1\over2}m_{b-1}}
\prod_{i=1}^{p'-3}
\left[{{1\over2}(m_{i-1}+m_{i+1}+\delta_{i,p'-a-1}+\delta_{i,b-1})
\atop m_i}\right]_q\!\!,
\end{array}
\end{displaymath}
where the sum is over all $\boldm=(m_1,m_2,\ldots,m_{p'-3})$
for which $\boldm\equiv\boldR(\boldu_{p'-a,b})$,
with $\tilde{\boldm}=(m_0,m_1,m_2,\ldots,m_{p'-3})$.
\end{theorem}

\Proof We proceed much as in the proof of Theorem \ref{FermThrm1}.
Let $t=p'-2$. Lemma \ref{InductLem} implies that:
\begin{displaymath}
\begin{array}{l}
\displaystyle
\chi^{1,p'+1-i}_{a_i,b_i,e_i,f_i}(m_{i-1},m_i)\\
\hskip8mm
\displaystyle
{}={}
\sum_{m\equiv m_{i-1}}
q^{{1\over4}\left((m_{i-1}-m_i)^2-\delta_{e_i+f_i,1}\right)}
\left[{{1\over2}(m_{i-1}+m)\atop m_i}\right]_q
\chi^{1,p'-i}_{a_{i+1},b_{i+1},e_i,f_i}(m_i,m)
\end{array}
\end{displaymath}
for $i=1,2,\ldots,t$. In the $i$th case, we replace the summation
variable $m$ with $m=m_{i+1}+\delta_{i,p'-a-1}+\delta_{i,b-1}$.
Thereupon, $m_{i-2}\equiv m_i+\delta_{i,p'-a}+\delta_{i,b}$
for $2\le i\le t+1$.
When $i\ne p'-a-1$ and $i\ne b-1$, we immediately obtain
$$
\chi^{1,p'-i}_{a_{i+1},b_{i+1},e_i,f_i}(m_{i},m_{i+1})
=\chi^{1,p'-i}_{a_{i+1},b_{i+1},e_{i+1},f_{i+1}}(m_{i},m_{i+1}),
$$
When $i=p'-a-1$, so that $a_{i+1}=p'-i-1$,
$e_i=0$ and $e_{i+1}=1$,
Lemma \ref{GenSwitchLem} yields:
$$
\chi^{1,p'-i}_{a_{i+1},b_{i+1},e_i,f_i}
(m_{i},m_{i+1}+1+\delta_{a+b,p'})
=\chi^{1,p'-i}_{a_{i+1},b_{i+1},e_{i+1},f_{i+1}+\delta_{a+b,p'}}
(m_{i},m_{i+1}+\delta_{a+b,p'}),
$$
and when $i=b-1$, so that $b_{i+1}=1$, $f_i=1$ and $f_{i+1}=0$,
Lemma \ref{GenSwitchLem} yields:
$$
\chi^{1,p'-i}_{a_{i+1},b_{i+1},e_i+\delta_{a+b,p'},f_i}(m_{i},m_{i+1}+1)
=q^{{1\over2}(a_b-1-m_{b-1})}
\chi^{1,p'-i}_{a_{i+1},b_{i+1},e_{i+1},f_{i+1}}(m_{i},m_{i+1}).
$$
By Lemma \ref{SeedLem},
$\chi^{1,2}_{1,1,1,0}(m_{t},m_{t+1})=\delta_{m_t,0}\delta_{m_{t+1},1}$,
so that we require $m_{t}=0$ and $m_{t+1}=1$
for a non-zero contribution,
whereupon the Gaussian polynomial in the $t$th summation is equal to 1.
Moreover, it then follows from 
$m_{i-2}\equiv m_i+\delta_{i,p'-a}+\delta_{i,b}$ that
$(m_1,m_2,m_3,\ldots,m_{t-1})\equiv\boldR(\boldu_{p'-a,b})$.
We now calculate
${1\over2}(a_b-1)-{1\over4}\sum_{i=1}^{t}\delta_{e_i+f_i,1}={1\over4}(a-b)$,
by using $e_i+f_i=1$ in $t-\vert t+2-a-b\vert$ cases and
$a_b=a$ if $b\le t-a+2$ and $a_b=t-b+2$ if $b\ge t-a+2$.
Combining all the above yields:
\begin{displaymath}
\begin{array}{l}
\displaystyle
\chi^{1,p'}_{a_1,b_1,e_1,f_1}(m_0,m_1)\\
\hskip2mm
\displaystyle
{}={}
\hskip-1mm
\sum
q^{{1\over4}\left(\sum_{i=1}^t(m_{i-1}-m_i)^2+a-b\right)-{1\over2}m_{b-1}}
\prod_{i=1}^{t-1}
\left[\!{{1\over2}(m_{i-1}\!+\!m_{i+1}\!+\!
\delta_{i,p'-a-1}\!+\!\delta_{i,b-1})
\atop m_i}\!\right]_q\!\!,
\end{array}
\end{displaymath}
where $m_t=0$ and the sum is over all
$(m_2,m_3,\ldots,m_{t-1})\!\equiv\!(R_2,R_3,\ldots,R_{t-1})$
when $\boldR(\boldu_{p'-a,b})=(R_1,R_2,\ldots,R_{t-1})$.
{}From $m_1\equiv R_1$, we obtain $m_0+e_1+f_1\equiv R_1$.
The theorem now follows after noting that
$$
\sum_{i=1}^t(m_{i-1}-m_i)^2=\tilde{\boldm}\boldC^{(t)}\tilde{\boldm}^T-L^2,
$$
and using Lemma \ref{FirstStepLem} in the form:
$$
\chi^{1,p'}_{a,b,c}(m_0)=\sum_{m_1\equiv R_1}
\chi^{1,p'}_{a_1,b_1,e_1,f_1}(m_0,m_1),
$$
and noting that $a=a_1$ and $b=b_1$.
\cqfd
\medskip

Again, we use Lemma \ref{ReflectLem} to obtain a further expression:

\begin{corollary}\label{FermCorol2}
Let $1\le a,b<p'$ and $c=b+1$.
Then, with $\boldC=\boldC^{(p'-2)}$, $m_0=L$, and $m_{p'-2}=0$,
\begin{displaymath}
\begin{array}{l}
\displaystyle
\chi^{1,p'}_{a,b,c}(L)\\
\hskip2mm
\displaystyle
{}={}
q^{{1\over4}(b-a-L^2)}
\hskip-1mm
\sum
q^{{1\over4}\tilde{\sboldm}\sboldC\tilde{\sboldm}^T-{1\over2}m_{p'-b-1}}
\hskip-2mm
\prod_{i=1}^{p'-3}
\left[{{1\over2}(m_{i-1}\!+\!m_{i+1}\!+\!\delta_{i,a-1}\!+\!\delta_{i,p'-b-1})
\atop m_i}\right]_q\!,\!
\end{array}
\end{displaymath}
where the sum is over all $\boldm=(m_1,m_2,\ldots,m_{p'-3})$
for which $\boldm\equiv\boldR(\boldu_{a,p'-b})$,
with $\tilde{\boldm}=(m_0,m_1,m_2,\ldots,m_{p'-3})$.
\end{corollary}

\section{ABF generating functions}\label{ABFFunSec}

We now use Lemma \ref{DualityLem} to convert the constant-sign expressions
obtained above to the ABF case.

\begin{theorem}
Let $1\le a,b<p'$,
$\boldC=\boldC^{(p'-3)}$, $m_0=L$ and $m_{p'-2}=0$.
Then if $b>1$ and $c=b-1$,
\begin{displaymath}
\begin{array}{l}
\displaystyle
\chi^{p'-1,p'}_{a,b,c}(L)\\
\hskip5mm
\displaystyle
{}={}
f_{a,b,c}
\sum_{\makebox[0pt]{$\scriptstyle\sboldm\equiv\sboldQ(\sboldu_{a,b})$}}
q^{{1\over4}{\sboldm}\sboldC{\sboldm}^T-{1\over2}m'_{a-1}}
\prod_{i=1}^{p'-3}
\left[{{1\over2}(m_{i-1}+m_{i+1}+\delta_{i,a-1}+\delta_{i,b-1})
\atop m_i}\right]_q\!,
\end{array}
\end{displaymath}
and also
\begin{displaymath}
\begin{array}{l}
\displaystyle
\chi^{p'-1,p'}_{a,b,c}(L)\\
\hskip5mm
\displaystyle
{}={}
f_{a,b,c}
\sum_{\makebox[0pt]{$\scriptstyle\sboldm\equiv\sboldR(\sboldu_{p'-a,b})$}}
q^{{1\over4}{\sboldm}\sboldC{\sboldm}^T-{1\over2}m'_{p'-a-1}}
\prod_{i=1}^{p'-3}
\left[{{1\over2}(m_{i-1}+m_{i+1}+\delta_{i,p'-a-1}+\delta_{i,b-1})
\atop m_i}\right]_q\!,
\end{array}
\end{displaymath}
and if $b<p'-1$ and $c=b+1$ then:
\begin{displaymath}
\begin{array}{l}
\displaystyle
\chi^{p'-1,p'}_{a,b,c}(L)\\
\hskip3mm
\displaystyle
{}={}
f_{a,b,c}
\sum_{\makebox[0pt]{$\scriptstyle\sboldm\equiv\sboldQ(\sboldu_{p'-a,p'-b})$}}
q^{{1\over4}{\sboldm}\sboldC{\sboldm}^T-{1\over2}m'_{p'-a-1}}
\hskip-2mm
\prod_{i=1}^{p'-3}
\left[{{1\over2}(m_{i-1}+m_{i+1}+\delta_{i,p'-a-1}+\delta_{i,p'-b-1})
\atop m_i}\right]_q\!,
\end{array}
\end{displaymath}
and also
\begin{displaymath}
\begin{array}{l}
\displaystyle
\chi^{p'-1,p'}_{a,b,c}(L)\\
\hskip5mm
\displaystyle
{}={}
f_{a,b,c}
\sum_{\makebox[0pt]{$\scriptstyle\sboldm\equiv\sboldR(\sboldu_{a,p'-b})$}}
q^{{1\over4}{\sboldm}\sboldC{\sboldm}^T-{1\over2}m'_{a-1}}
\prod_{i=1}^{p'-3}
\left[{{1\over2}(m_{i-1}+m_{i+1}+\delta_{i,a-1}+\delta_{i,p'-b-1})
\atop m_i}\right]_q\!,
\end{array}
\end{displaymath}
where $f_{a,b,c}=q^{-{1\over4}(a-b)(a-c)}$, and
in each case $m'_i=m_i$ for $i>0$ and $m'_0=0$.
\end{theorem}

\Proof Using
$\left[{P\atop Q}\right]_{q^{-1}}
=q^{-Q(P-Q)}\left[{P\atop Q}\right]_{q}$ yields:
\begin{displaymath}
\begin{array}{l}
\displaystyle
\prod_{i=1}^{p'-3}
\left[{{1\over2}(m_{i-1}+m_{i+1}+\delta_{i,a-1}+\delta_{i,b-1})
\atop m_i}\right]_{q^{-1}}\\
\hskip2mm
\displaystyle
{}={}
\hskip-1.5mm
\prod_{i=1}^{p'-3}
q^{-{1\over2}m_i(m_{i\!-\!1}-2m_i+m_{i\!+\!1}
       +\delta_{i,a\!-\!1}+\delta_{i,b-\!1})}
\hskip-1mm
\left[{{1\over2}(m_{i-1}+m_{i+1}+\delta_{i,a-1}+\delta_{i,b-1})
\atop m_i}\right]_{q}\\
\hskip2mm
\displaystyle
{}={}
q^{{1\over2}({\sboldm}\sboldC{\sboldm}^T-m_1m_0-m'_{a-1}-m'_{b-1})}
\prod_{i=1}^{p'-3}
\left[{{1\over2}(m_{i-1}+m_{i+1}+\delta_{i,a-1}+\delta_{i,b-1})
\atop m_i}\right]_{q}\!.
\end{array}
\end{displaymath}
Substituting this and
$
\tilde{\boldm}\boldC^{(p'-2)}\tilde{\boldm}^T=
{\boldm}\boldC{\boldm}^T
+2L^2-2m_0m_1
$
into Theorem \ref{FermThrm1}, and noting that $b>1$,
we obtain:
\begin{displaymath}
\begin{array}{l}
\displaystyle
\chi^{1,p'}_{a,b,c}(L;q^{-1})\\
\hskip5mm
\displaystyle
{}={}
q^{{1\over4}(b-a-L^2)}
\sum
q^{{1\over4}{\sboldm}\sboldC{\sboldm}^T-{1\over2}m'_{a-1}}
\prod_{i=1}^{p'-3}
\left[{{1\over2}(m_{i-1}+m_{i+1}+\delta_{i,a-1}+\delta_{i,b-1})
\atop m_i}\right]_q\!.
\end{array}
\end{displaymath}
The first expression then follows from Lemma \ref{DualityLem}.
The other three expressions arise in a similar way from
Theorem \ref{FermThrm2}, Corollary \ref{FermCorol1} and
Corollary \ref{FermCorol2} respectively.
\cqfd

\begin{appendix}

\section{$\boldm\boldn$-systems}

Consider once more the proof of Lemma \ref{InductLem}.
There, It is shown that for each path $h'\in{\P}^{p'}_{a,b,e,f}(L,m)$,
there is a unique pair $(h,k)$ with
$h\in{\P}^{p'-1}_{a-e,b-f,e,f}(m',m)$ and $k\in\N$, for which $h'$
arises from the action of a ${\B}$-transform on $h$,
followed by the insertion of $k$ particles,
followed by moving these particles in some way.
The path $h$ will be referred to as the {\it $(e,f)$-antecedent} of $h'$.
The value of $k$ will be referred to as the
{\it $(e,f)$-particle content} of $h'$.
{}From the proof of Lemma \ref{InductLem}, we see that it is given by:
\begin{equation}\label{OneContentEq}
2k=L'-2m'+m.
\end{equation}

Now, given $e_i,f_i\in\{0,1\}$ for $1\le i\le t=p'-2$,
we may iterate the above procedure.
Let $h\in{\P}^{p'}_{a,b,e,f}(L,m)$, and let $n_1$ be its
$(e_1,f_1)$-particle content and $h'$ its $(e_1,f_1)$-antecedent
(we deviate from the above priming convention).
Now let $n_2$ and $h''$ be respectively the $(e_2,f_2)$-particle content
and $(e_2,f_2)$-antecedent of $h'$.
Proceeding in this way, we obtain a vector
$\boldn=(n_1,n_2,\ldots,n_{t})$,
which we shall simply refer to as the particle content of $h$.
Note that the particle content is dependent on the
particular sequence of $e_i,f_i\in\{0,1\}$ being considered.
Thus, in general, a given path $h$ has differing particle contents
in the two systems considered in Sections \ref{FirstSystemSec} and
\ref{SecondSystemSec}.


\subsection{First system}\label{AFirstSystemSec}

On examining the proof of Theorem \ref{FermThrm1}, we see that the
use of the $i$th ${\B}$-transform therein, results from
substituting $L'=m_{i-1}$, $m'=m_i$ and
$m=m_{i+1}+\delta_{i,a-1}+\delta_{i,b-1}$ into Lemma \ref{InductLem}.
Therefore, with $e_i$ and $f_i$ given by (\ref{EFEq1}), we find that
if a path $h\in{\P}^{p'}_{a,b,c}(L)$ has particle content
$(n_1,n_2,\ldots,n_t)$ then, from (\ref{OneContentEq}),
\begin{equation}\label{MNEq1}
m_{j-1}+m_{j+1}=2m_j+2n_j-\delta_{j,a-1}-\delta_{j,b-1},
\end{equation}
for $1\le j\le t$, where $m_0=L$ and $m_t=m_{t+1}=0$.

The set of $t$ equations given by (\ref{MNEq1})
defines an interdependence between the vectors
$\boldn=(n_1,n_2,\ldots,n_t)$ and
$\boldm=(m_0,m_1,\ldots,m_{t-1})$
known as the $\boldm\boldn$-system.

If we define the vector $\boldu_{a,b}=(u_1,u_2,\ldots,u_t)$
with components $u_j=\delta_{j,a-1}+\delta_{j,b-1}$,
the $\boldm\boldn$-system described in this section may
be conveniently written:
\begin{equation}
2\boldn=-\boldm\boldC^{(t)}+\boldu_{a,b},
\end{equation}
where $\boldC^{(t)}$ is the Cartan matrix of type $A_t$,
defined by (\ref{CartanDef}).

\subsection{Second system}\label{ASecondSystemSec}

By the same means as in Section \ref{AFirstSystemSec},
we obtain the $\boldm\boldn$-system for the case considered
in Section \ref{SecondSystemSec}.
Thus, with $e_i,f_i$ defined by (\ref{EFEq2}),
the proof of Theorem \ref{FermThrm2} shows that if a path
$h\in{\P}^{p'}_{a,b,c}(L)$ has particle content
$(n_1,n_2,\ldots,n_t)$ then
\begin{equation}\label{MNEq2}
m_{j-1}+m_{j+1}=2m_j+2n_j-\delta_{j,p'-a-1}-\delta_{j,b-1}-\delta_{j,t},
\end{equation}
for $1\le j\le t$, where $m_0=L$ and $m_t=m_{t+1}=0$.
We then obtain:
\begin{equation}
2\boldn=-\boldm\boldC^{(t)}+\boldu_{p'-a,b}+\bolde_t,
\end{equation}
where the $t$-dimensional $\bolde_t=(0,0,\ldots,0,1)$.

\section{Bijection between paths and partitions}

In this appendix, we briefly describe a natural weight-preserving
bijection between the paths that we've been considering in this
paper and partitions that satisfy certain hook-difference
conditions \cite{abbbfv,flpw}.
Such a bijection occurs in both the parafermion and the ABF cases.
In fact, as pointed out in \cite{abbbfv}, these bijections are
just special cases of a more general bijection existing
for the weighted paths of \cite{forrester-baxter}.
We give a full description (of the general case) in \cite{flpw}.

\subsection{Partitions with prescribed hook-difference conditions}

A partition $\mu=(\mu_1,$$\mu_2,$$\ldots,$$\mu_M)$ is a sequence 
of $M$ integer parts $\mu_1,\mu_2,\ldots,\mu_M$ satisfying
$\mu_1\ge\mu_2\ge\cdots\ge\mu_M\ge0$.
The weight $\wt(\mu)$ of $\mu$ is given by $\wt(\mu)=\sum_{i=1}^M \mu_i$.
The partition $\mu$ is often depicted by its {\it Young diagram}
(also called Ferrars graph),
$F^\mu$ which comprises $M$ left-adjusted
rows, the $i$th row of which (reading down) consists of $\mu_i$
cells \cite{andrews-red-book}.
The coordinate $(i,j)$ of a cell is obtained by setting $i$
and $j$ to be respectively, the row and column (reading from the left)
in which the cell resides.
The $k$th diagonal of $F^\mu$ comprises all those cells of $F^\mu$
with coordinates $(i,j)$ which satisfy $i-j=k$.

The partition $\mu^\prime$, conjugate to $\mu$, is obtained
by setting $\mu^\prime_j$ to be the number of cells in the
$j$th column of $F^\mu$.
The hook-difference at the cell with coordinate $(i,j)$
is then defined to be $\mu_i-\mu^\prime_j$.
As an example, filling each cell of $F^{(5,4,3,1)}$ with its hook difference,
yields:
$$
\youngd{
\multispan{11}\hrulefill\cr
&1&&\mathbf2&&2&&3&&4&\cr
\multispan{11}\hrulefill\cr
&0&&1&&\mathbf1&&2&\cr
\multispan{9}\hrulefill\cr
&-1&&0&&0&\cr
\multispan{7}\hrulefill\cr
&-3&\cr
\multispan{3}\hrulefill\cr}\:.
$$

\medskip\noindent
The bold entries are those on diagonal $-1$.
In what follows, we will be especially interested in the
hook-differences on certain diagonals.

Let $K,i,N,M,\alpha,\beta$ be non-negative integers for
which $1\le i\le K/2$, $\alpha+\beta<K$ and
$\beta-i\le N-M\le K-\alpha-i$.
In \cite{abbbfv}, $D_{K,i}(N,M;\alpha,\beta)$ is defined to be
the generating function for partitions $\mu$
into at most $M$ parts, each not exceeding $N$ such that
the hook differences on diagonal $1-\beta$ are at least
$\beta-i+1$, and on diagonal $\alpha-1$ are at most $K-i-\alpha-1$.
In addition, if $\alpha=0$, the restriction that $\mu_{N-L+i+1}>0$
is also imposed; and if $\beta=0$, the restriction that
$\mu_1>M-i$ is also imposed.

\subsection{Parafermions}

In the parafermion case,
the bijective image of a path $h\in{\P}^{p'}_{a,b,c}(L)$ is
obtained as follows.
Begin with an empty Young diagram.
Now traverse the path from left to right.
If a scoring vertex that contributes $x$ to the weight is
encountered, append a new first row of length $x$ to the
top of the Young diagram.
If a scoring vertex that contributes $y$ to the weight is
encountered, append a new first column of length $y$ to the
left of the Young diagram.
The diagram is not changed at non-scoring vertices.
In this way, after all $L$ vertices have been considered,
a Young diagram $F^\mu$ results.

For example, consider the path shown in Fig.~1, and let $c=b+1$.
Here, we obtain the following Young diagram:

\medskip
\centerline{\epsfig{file=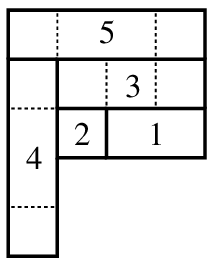}\:\raisebox{30pt}[0pt].}
\medskip

\noindent Here, the entries indicate the order in which the
pieces have been added.

If $c=b+1$ let $r=0$, and if $c=b-1$ let $r=1$. It may be shown 
(\cite{flpw}, or by using the techniques of \cite{foda-warnaar}) 
that each hook-difference on diagonal $1-r$ of $\mu$ is at least 
$r-a+1$, and each hook-difference on diagonal $-r$ of $\mu$ is at 
most $r-a+p'$. Moreover, it may be shown that the map is in fact 
a bijection. Since clearly $\wt(\mu)=\wt(h)$, we have
$$
\chi^{1,p'}_{a,b,c}(L)=D_{p',a}\left({L-a+b\over2},{L+a-b\over2},1-r,r\right).
$$

\subsection{ABF}

The description of the bijection in the ABF case proceeds similarly.
We must first provide an analogue of $c$ that was defined
in Section \ref{AltPresSec}.
If the $i$th vertex has coordinates $(x,y)$,
define $\tilde c(h_{i-1},h_i,h_{i+1})$ as follows:
\begin{eqnarray*}
\tilde c(h-1,h,h+1)&=&x\,;\\
\tilde c(h+1,h,h-1)&=&y\,;\\
\tilde c(h-1,h,h-1)&=&0\,;\\
\tilde c(h+1,h,h+1)&=&0\,.
\end{eqnarray*}
Now, if we define
\begin{equation}\label{WtIIIDef}
\wt^{\rm III}(h)=\sum_{i=1}^L \tilde c_i(h_{i-1},h_i,h_{i+1}),
\end{equation}
then, as is readily shown,
\begin{equation}
\chi^{p'-1,p'}_{a,b,c}(L)
=\sum_{h\in{\P}^{p'}_{a,b,c}(L)} q^{\wt^{\rm III}(h)},
\end{equation}
where $\chi^{p'-1,p'}_{a,b,c}(L)$ is given by (\ref{RenormABFEq})
and (\ref{GenIIIDef}).

The partition $\mu$ is now obtained exactly as in the above
description of the
parafermion case:
the only difference being that the non-scoring vertices
there are scoring vertices here, and vice-versa.

For example, again consider the path shown in Fig.~1, and let $c=b+1$.
In the ABF case, we obtain the following Young diagram:

\medskip
\centerline{\epsfig{file=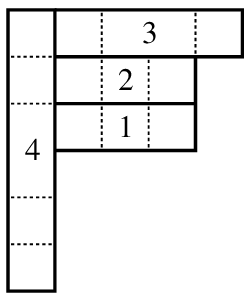}\:\raisebox{30pt}[0pt].}
\medskip

Let $r=\min(b,c)$.
It may be shown that each hook-difference on diagonal
$1-r$ of $\mu$ is at least $r-a+1$, and each hook-difference on diagonal
$p'-2-r$ of $\mu$ is at most $r-a+p'$.
Once more, the map may be shown to be a bijection, whereupon:
$$
\chi^{p'-1,p'}_{a,b,c}(L)
=D_{p',a}\left({L-a+b\over2},{L+a-b\over2},p'-r-1,r\right).
$$
(cf.~eq.~(5.1) of \cite{abbbfv}).

Finally, we note that this bijection generalises that given 
in \cite{foda-warnaar}, and that the method given there may 
be readily extended to deal with the current case.

\end{appendix}

\subsection*{Acknowledgments}
 
We wish to thank Keith Lee and Slava Pugai for collaboration on 
\cite{flpw}, on which this work is directly based, and for many 
stimulating remarks and discussions. We thank Ole Warnaar for
bringing \cite{schilling} to our attention. One of us (OF) also wishes 
to thank Naihuan Jing and Kailash Misra for the invitation to 
present this work in {Affine Lie Algebras and related topics},
North Carolina State University, Raleigh, 1998, and for their 
excellent hospitality. This work was supported by the Australian 
Research Council.

\end{document}